\documentclass[11 pt,a4paper]{article}

\usepackage[utf8]{inputenc}
\usepackage[english]{babel}   
\usepackage{graphicx}
\usepackage{lmodern}
\usepackage{amsmath}
\usepackage{tikz}
\usetikzlibrary{calc}
\usepackage{subfig}
\usepackage{longtable}
\usepackage{epsfig}
\usepackage{color}
\usepackage{setspace}
\usepackage[top=2cm,bottom=2.5cm,left=2.5cm,right=2.5cm]{geometry}
\usepackage{sectsty}
\usepackage{fancyhdr}
\usepackage{wasysym}
\usepackage[colorlinks=true]{hyperref}
\usepackage{caption}
\usetikzlibrary{patterns}
\usepackage{amsfonts}
\usepackage{amssymb}
\usepackage{amsthm}
\usepackage{multicol}

\definecolor{vert}{rgb}{0.2,0.6,0}
                 
\newtheorem{theorem}{Theorem}
\newtheorem{remarque}[theorem]{Remark}

\newcommand{\Tr}{^{\mathbf{T}}}

\newcommand{\veps}{\varepsilon}

\title{Quasi-static crack propagation with a {G}riffith criterion using a variational discrete element method}
\author{\begin{minipage}{\textwidth}\centering Fr\'ed\'eric
Marazzato$^{1,2,3}$, Alexandre Ern$^{2,3}$ and Laurent Monasse$^{4}$\\
   \small{$^{1}$Department of Mathematics, Louisiana State University, Baton Rouge, LA 70803, USA}\\
   \small{email: \texttt{marazzato@lsu.edu}}\\
   \small{$^{2}$CERMICS, Ecole des Ponts, 77455 Marne-la-Vall\'ee, France}\\
   \small{email: \texttt{alexandre.ern@enpc.fr}}\\
   \small{$^{3}$Inria, 2 rue Simone Iff, 75589 Paris, France}\\
   \small{$^4$Universit\'e C\^ote d'Azur, Inria, CNRS, LJAD, EPC COFFEE, 06108 Nice, France}\\
   \small{email: \texttt{laurent.monasse@inria.fr}}\end{minipage}}
   
\date{}
   
\begin{document}
\hypersetup{urlcolor=blue,linkcolor=red,citecolor=blue}

\maketitle

\begin{abstract}
A variational discrete element method is applied to simulate quasi-static crack propagation.
Cracks are considered to propagate between the mesh cells through the mesh facets. The elastic behaviour is parametrized by the continuous mechanical parameters (Young modulus and Poisson ratio).
A discrete energetic cracking criterion coupled to a discrete kinking criterion guide the cracking process.
Two-dimensional numerical examples are presented to illustrate the robustness and versatility of the method.
\end{abstract}

\section{Introduction}
Discrete element methods (DEM) are popular in the modeling of granular materials, soil and rock mechanics. 
DEM generally use sphere packing to discretize the domain as small spheres interacting through forces and torques \cite{MR2548422}, but the main difficulty is to derive a suitable set of parameter values for those interactions so as to reproduce a given Young modulus $E$ and Poisson ratio $\nu$ at the macroscopic level \cite{jebahi2015discrete,MR3735741}.
Advantages of DEM are their ability to deal with discontinuous materials, such as fractured or porous materials, as well as the
possibility to take advantage of GPU computations \cite{MR3606237}.
A first DEM parametrized only by $E$ and $\nu$ has been proposed in \cite{LM_CM_2012} for 
elastic computations on Voronoi meshes.
In a consecutive work \cite{MaraDEM}, a variational DEM has been proposed for
elasto-plasticity computations on polyhedral meshes using cell-wise reconstructions of the strains.
The numerical results reported in \cite{MaraDEM} confirmed that the
macroscopic behaviour of elastic continua is indeed correctly reproduced by the 
variational DEM. The method developed in
\cite{MaraDEM} takes its roots in \cite{eymard2009discretization} which is indeed a hybrid finite volume method. It is called variational DEM since it is possible to reinterpret the method as a consistent discretization of elasto-plasticity with discrete elements. In particular, a force-displacement interpretation of the method is derived from the usual stress-strain approach. Also, the mass matrix is diagonal and the stencil for the gradient reconstruction is compact as in usual DEM.

DEM for cracking have been developed in \cite{andre2013using} and \cite{andre2019novel} with cracks propagating through the facets of the (Voronoi) mesh and using a critical stress criterion (initiation criterion). 
Coupled FEM-DEM techniques for crack computations, as \cite{zarate2015simple} (2d) and \cite{zarate2018three} (3d), have been introduced to take advantage of the FEM ability in computing elasticity and of the ability of DEM to handle cracked media.
A similar approach, but using a different reconstruction of strains based on moving least-squares interpolations, can be traced back to \cite{belytschko1994element} (2d) and \cite{sukumar1997element} (3d). 
Crack propagation can be based instead on the Griffith criterion which relies on the computation of the stress intensity factors (SIF) at the crack tip when coupled with the Irwin formula. 
Virtual element methods (VEM) have been recently applied to crack propagation \cite{hussein2018virtual}. Cracks were allowed to cut through the polyhedral mesh cells as in the extended finite element method (XFEM) which is based on an extended space of basis functions \cite{chahine2008crack} and a level-set description of the crack \cite{moes2002x}. Phase-field methods instead smooth the crack and have been developed among others in \cite{bourdin2000numerical} and subsequent work.
Phase-field methods are not based on SIF computations but rather on a variational formulation of cracking \cite{francfort1998revisiting}. 
Furthermore, DEM using cohesive laws have been developed for fragmentation computations \cite{mariotti2009modeling} with a view towards uniting initiation and propagation. These methods allow one to devise an initiation criterion and also to control the energy dissipation as with a Griffith criterion. The cracks still go through the mesh facets. This is also the case for similar methods of higher-order such as discontinuous Galerkin methods \cite{hansbo2015discontinuous}.

The main goal of the present work is to develop a variational DEM using a Griffith criterion to compute crack propagation through the mesh facets. The method supports in principle polyhedral meshes, but the present numerical experiments are restricted to triangular meshes.
The proposed method is close to \cite{MaraDEM} (where there is no cracking) but the degrees of freedom (dofs) are different.
Only cell dofs are used in the present work.
The cracking algorithm hinges on two main ingredients. The first ingredient is an approximation of the energy release rate at every vertex along the crack.
The second ingredient is a kinking criterion used to determine the next breaking facet and thus the crack path. The kinking criterion, in the spirit of \cite{sih1974strain}, consists in selecting for the crack path the inner facet of the mesh that maximizes a quantity representing the local density of elastic energy.

The present work is organized as follows. Section~\ref{sec:governing cracking} briefly recalls the equations of elasticity and cracking in a Cauchy continuum. Section~\ref{sec:space discretization} 
introduces the proposed variational DEM and presents the space discretization of the governing equations.
Moreover, a numerical test is reported to assess the convergence of the space discretization in the presence of a singularity.
Section~\ref{sec:quasi-static cracking} addresses the full discretization of the quasi-static cracking problem. Section~\ref{sec:numerics} contains numerical results on quasi-static crack propagation problems in two space dimensions. 
Finally, Section~\ref{sec:conclusions cracking} draws some conclusions.

\section{Governing equations for quasi-static cracking}
\label{sec:governing cracking}
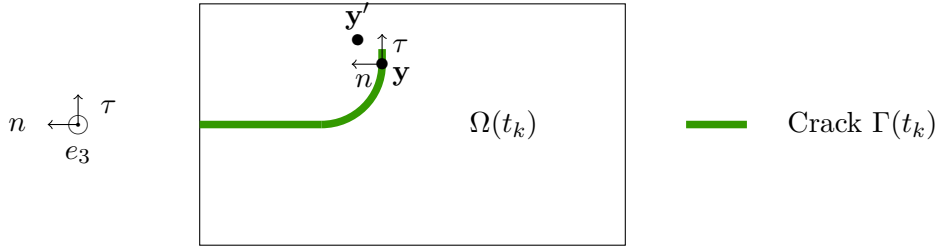
\begin{figure} [!htp]
\begin{center}
\begin{tikzpicture} [scale = 0.8]
\pgfmathsetmacro{\R}{1}
\pgfmathsetmacro{\length}{7}
\pgfmathsetmacro{\width}{2}
\pgfmathsetmacro{\lz}{2}
\draw (\length,\width) -- (\length,-\width) -- (0,-\width) -- (0,\width) -- cycle;
\draw[vert,line width=1mm] (0,0) -- (\lz,0);
\draw[vert,line width=1mm] (\lz,0) arc (-90:0:\R);
\draw[vert,line width=1mm] (\lz+\R,\R) -- (\lz+\R,\R+0.25);

\draw (\lz+\R,\R) node {$\bullet$};
\draw (\lz+\R,\R-0.2) node[right] {$\mathbf{y}$};
\draw[->] (\lz+\R,\R) -- (\lz+\R-0.5,\R);
\draw (\lz+\R-0.3,\R) node[below] {$n$};
\draw[->] (\lz+\R,\R) -- (\lz+\R,\R+0.5);
\draw (\lz+\R,\R+0.3) node[right] {$\tau$};
\draw (\lz+\R-0.4,\R+0.8) node {$\mathbf{y}'$};
\draw (\lz+\R-0.4,\R+0.4) node {$\bullet$};

\draw[->] (-2,0) -- (-2.5,0);
\draw (-3,0) node {$n$};
\draw[->] (-2,0) -- (-2,0.5);
\draw (-1.5,0.3) node {$\tau$};
\draw (-2,0) node {$\odot$};
\draw (-2.,-0.5) node{$e_3$};

\draw[vert,line width=1mm] (8, 0) -- (9, 0);
\draw[right] (9.5, 0) node {Crack $\Gamma(t_k)$};
\draw (5,0) node {$\Omega(t_k)$};
\end{tikzpicture}
\caption{Sketch of a crack in the two-dimensional domain $\Omega(t_k)$.}
\label{fig:3d crack}
\end{center}
\end{figure}

We consider an elastic fragile material occupying the domain
$\Omega \subset \mathbb{R}^2$ in the reference configuration
and evolving over the finite pseudo-time interval $[0,T]$, $T > 0$, under the action of a volumetric force $f$ and boundary conditions. The pseudo-time interval $[0,T]$ is discretized by means of $(K+1)$ 
discrete pseudo-time nodes $(t_k)_{k\in\{0,\ldots,K\}}$ with $t_0:=0$ and $t_{K}:=T$. 
The strain regime is restricted to small strains so that we use the linearized strain tensor $\veps(u):=\frac12(\nabla u+(\nabla u)\Tr) \in \mathbb{R}^{2\times 2}$, where $u$ is the $\mathbb{R}^2$--valued displacement field.
The material is supposed to be homogeneous and isotropic.  
The stress tensor $\sigma(u) \in \mathbb{R}^{2\times 2}$ is such that
\begin{equation}
\sigma(u) := \mathbb{C}: \veps(u),
\end{equation}
where $\mathbb{C}$ is the fourth-order stiffness tensor. The elastic material is characterized by the Young modulus $E$ and the Poisson ratio $\nu$ or equivalently by the Lam\'e coefficients $\lambda$ and $\mu$. 
The boundary of $\Omega$ is partitioned as $\partial \Omega = \partial \Omega_D \cup \partial \Omega_N$, a Dirichlet condition is prescribed on $\partial\Omega_D$, and a Neumann condition on $\partial \Omega_N$, so that we enforce for all $k=0,\cdots,K$,
\begin{equation}
u = u_D(t_k) \ \text{ on } \partial \Omega_D, \qquad
\sigma(u) \cdot n = g_N(t_k)\  \text{ on } \partial \Omega_N. 
\end{equation}
Since cracking can occur, we denote $\Gamma(t_k)$ the crack at the pseudo-time node $t_k$ and the actual domain at the pseudo-time node $t_k$ is 
\begin{equation}
\Omega(t_k) := \Omega \setminus \Gamma(t_k).
\end{equation} 
This implies that $\partial \Omega(t_k) = \partial \Omega_D \cup \partial \Omega_N \cup \Gamma(t_k)$. We enforce a homogeneous Neumann condition on $\Gamma(t_k)$ for all $k=0,\cdots,K$, i.e.,
\begin{equation}
\sigma(u) \cdot n = 0\  \text{ on } \Gamma(t_k). 
\end{equation}
Since we are interested in crack propagation, we assume that $\Omega(0)$ already contains a crack, i.e., $\Gamma(0)\ne\emptyset$. 
The crack $\Gamma(t_k)$ is supposed to be a countably rectifiable
1--manifold for all $k=0,\cdots,K$ (see \cite{dal2013generalised}). 
This hypothesis ensures the almost everywhere (a.e.) existence of a normal vector $n$ and a tangent vector $\tau$ to $\Gamma(t_k)$ at any point $\mathbf{y} \in \Gamma(t_k)$ \cite{simon1983lectures}.
Figure \ref{fig:3d crack} illustrates these quantities.

The stress intensity factors (SIF) at any point $\mathbf{y} \in \Gamma(t_k)$
are usually defined for a purely elastic material as
\begin{equation}
\left\{
\begin{aligned}
& K_1(\mathbf{y}) := \mathop{\mathrm{lim}} \limits_{\mathbf{y}' \to \mathbf{y}} \sigma_{nn}(\mathbf{y}') \sqrt{2 \pi d(\mathbf{y},\mathbf{y}')}, \\
& K_2(\mathbf{y}) := \mathop{\mathrm{lim}} \limits_{\mathbf{y}' \to \mathbf{y}} \sigma_{n\tau}(\mathbf{y}') \sqrt{2 \pi d(\mathbf{y},\mathbf{y}')}, \\
\end{aligned}
\right.
\label{eq:FIC contrainte 3d}
\end{equation}
where $d(\cdot,\cdot)$ is the Euclidean distance in $\mathbb{R}^2$.
If the stresses remain bounded in the vicinity of $\mathbf{y} \in \Gamma(t_k)$, then the SIF are null.
Using the Irwin formula, one can define the energy release rate $\mathcal{G}(\mathbf{y})$ in the plane strain hypothesis as
\begin{equation}
\label{eq:Irwin formula}
\mathcal{G}(\mathbf{y}) := \frac{1-\nu^2}{E}\left(K_1(\mathbf{y})^2 + K_2(\mathbf{y})^2 \right). 
\end{equation}
Admissible states are characterized by the inequality
\begin{equation}
\mathcal{G}(\mathbf{y}) \leq \mathcal{G}_c, \quad \forall \mathbf{y} \in \Gamma(t_k),
\label{eq:continuous fracture criterion}
\end{equation}
where $\mathcal{G}_c$ is a material property associated with the capacity of the material to sustain loads without locally failing and thus opening cracks. 
The material remains healthy at the point $\mathbf{y} \in \Gamma(t_k)$
if $\mathcal{G}(\mathbf{y}) < \mathcal{G}_c$
and breaks if $\mathcal{G}(\mathbf{y}) = \mathcal{G}_c$. The material parameter 
$\mathcal{G}_c$ is assumed to be homogeneous for simplicity.

To formulate the governing equations for quasi-static cracking,
we consider the following functional spaces depending on the pseudo-time node $t_k$: 
\begin{equation}
V_D(t_k) := \left\{ v \in H^1(\Omega(t_k);\mathbb{R}^d) \ | \ v_{|\partial \Omega_D} = u_D(t_k) \right\},
\qquad
V_0(t_k) := \left\{ v \in H^1(\Omega(t_k);\mathbb{R}^d) \ | \ v_{|\partial \Omega_D} = 0 \right\},
\end{equation}
where standard notation is used for the Hilbert Sobolev spaces. 
The weak solution is searched as a pair $(u,\Gamma)$ such that for all $k=0,\cdots,K$, $u(t_k) \in V_D(t_k)$, $\Gamma(t_k)\subset \Omega$ is a
1--manifold satisfying the above assumptions, and 
\begin{equation}
\label{eq:cracking-weak}
\left\{ \begin{alignedat}{2}
&a(t_k;u(t_k),\tilde{v}) = l(t_k;\tilde{v}),
&\quad&\forall \tilde{v}\in V_0(t_k),\\
&\mathcal{G}(\mathbf{y}) \leq \mathcal{G}_c, &\quad&\forall \mathbf{y} \in \Gamma(t_k).
\end{alignedat} \right.
\end{equation}
Here we introduced the stiffness bilinear form such that for all $(v,\tilde{v}) \in V_D(t_k)\times V_0(t_k)$, 
\begin{equation}
\label{eq:bilinear form crack}
a(t_k;v,\tilde{v}) := \int_{\Omega(t_k)} \varepsilon(v) : \mathbb{C} : \varepsilon(\tilde{v}), 
\end{equation}
and the linear form acting on $V_0(t_k)$ as follows:
\begin{equation}
\label{eq:linear form crack}
l(t_k;\tilde{v}) := \int_{\Omega(t_k)} f(t_k) \cdot \tilde{v} + \int_{\partial \Omega_N} g_N(t_k) \cdot \tilde{v}. 
\end{equation}
Note that the Dirichlet condition on $\partial\Omega_D$ is enforced strongly, whereas the Neumann condition on $\partial\Omega_N \cup \Gamma(t_k)$ is enforced weakly. 

\section{Space semi-discretization}
\label{sec:space discretization}
In this section, we present the space semi-discretization of~\eqref{eq:cracking-weak} using a variational DEM. 

\subsection{Discrete sets and degrees of freedom}

The domain $\Omega$ is discretized with a mesh $\mathcal{T}_h$ of size $h$ made of polygons with straight edges. 
We assume that $\Omega$ is itself a polygon so that the mesh covers $\Omega$ exactly. 
We also assume that the mesh
is compatible with the initial crack position $\Gamma(0)$ and
with the partition of the boundary into the Dirichlet and Neumann parts. 
Recall that the space dimension is $d=2$.

Let $\mathcal{C}$ denote the set composed of the mesh cells and, for all $k=0,\dots,K$, let $\mathcal{F}(t_k)$ denote the set composed of the mesh facets. This set depends on the pseudo-time node $t_k$ since a facet $F \in \mathcal{F}(t_k)$ is replaced, after cracking, by two boundary facets $F_-,F_+ \in \mathcal{F}(t_k)$ 
($F_-,F_+$ are the same geometric object, but are different objects regarding the data structure since each one belongs to the boundary of a different mesh cell). The barycentre of a mesh cell $c \in \mathcal{C}$ is denoted by $\mathbf{x}_c$ and the barycentre of a mesh facet $F \in \mathcal{F}(t_k)$ is denoted by $\mathbf{x}_F$.

Let $t_k$ be a pseudo-time node with $k=0,\cdots,K$.
We partition the set of mesh facets as $\mathcal{F}(t_k)=\mathcal{F}^i(t_k)\cup \mathcal{F}^b(t_k)$, where 
$\mathcal{F}^i(t_k)$ is composed of the internal facets shared
by two mesh cells and $\mathcal{F}^b(t_k)$ is the collection of the boundary facets 
sitting on the boundary $\partial\Omega(t_k)= \partial \Omega_D \cup \partial \Omega_N 
\cup \Gamma_{h}(t_k)$, where $\Gamma_h(t_k)$ denotes the discrete crack at $t_k$. 
Notice that every boundary facet belongs to the boundary of only one mesh cell. 
The subsets $\mathcal{F}^i(t_k)$ and $\mathcal{F}^b(t_k)$ depend on the pseudo-time node $t_k$ since, as the facet 
$F \in \mathcal{F}^i(t_k)$ cracks, it is replaced by the facets $F_+,F_- \in \mathcal{F}^b(t_k)$.
The discrete crack $\Gamma_h(t_k)$ is 
composed of facets belonging to a subset of $\mathcal{F}^b(t_k)$.
This subset is denoted by $\mathcal{F}^\Gamma(t_k) \subset \mathcal{F}^b(t_k)$. We also introduce the partition between boundary facets with Neumann boundary conditions $\mathcal{F}^b_N(t_k)$ (recall that homogeneous Neumann boundary conditions are imposed on newly created crack lips) and with Dirichlet boundary conditions $\mathcal{F}^b_D$ which does not depend on $t_k$. One thus has $\mathcal{F}^b(t_k) = \mathcal{F}^b_N(t_k) \cup \mathcal{F}^b_D$.

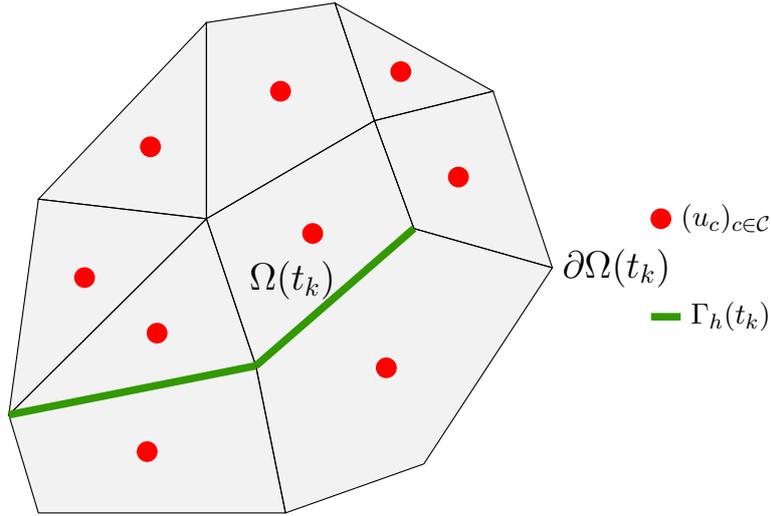
\begin{figure} [!htp] 
\begin{center}
\begin{tikzpicture}[scale=1.3]
\coordinate (a) at (0,-1);
\coordinate (b) at (0,1);
\coordinate (c) at (1.3,1.2);
\coordinate (d) at (-1.7,-0.8);
\coordinate (e) at (-2,-3);
\coordinate (f) at (0.5,-2.5);
\coordinate (g) at (1.7,0);
\coordinate (h) at (2.9,0.3);
\coordinate (i) at (2.1,-1.1);
\coordinate (j) at (3.5,-1.5);
\coordinate (k) at (2.2,-3.5);
\coordinate (l) at (0.8,-4);
\coordinate (m) at (-1.7,-4);

\draw[fill=gray,opacity=0.1] (b) --(c) -- (h) -- (j) -- (k) -- (l) --
(m) -- (e) -- (d) -- cycle;
\node[] at (barycentric cs:f=0.8,a=0.25,g=0.25,i=0.25) {\Large
  $\Omega(t_k)$};
\node[right] at (j) {\Large $\partial\Omega(t_k)$};
\path[draw] (a)-- (b)-- (c) -- (g) -- cycle;
\path[draw] (a)-- (b)-- (d) -- cycle;
\path[draw] (a) -- (e) -- (d) -- cycle;
\path[draw] (a) -- (f) -- (e) -- cycle;
\path[draw] (a) -- (g) -- (i) -- (f) -- cycle;
\path[draw] (c) -- (h) -- (g) -- cycle;
\path[draw] (g) -- (h) -- (j) -- (i) -- cycle;
\path[draw] (i) -- (j) -- (k) -- (l) -- (f) -- cycle;
\path[draw] (f) -- (l) -- (m) -- (e) -- cycle;

\path[draw,line width=1mm,vert] (i) -- (f);
\path[draw,line width=1mm,vert] (e) -- (f);

\fill[red] (barycentric cs:a=0.3,b=0.3,d=0.3) circle (3pt);
\fill[red] (barycentric cs:a=0.3,e=0.3,d=0.3) circle (3pt);
\fill[red] (barycentric cs:f=0.3,e=0.3,a=0.3) circle (3pt);
\fill[red] (barycentric cs:a=0.25,g=0.25,c=0.25,b=0.25) circle (3pt);
\fill[red] (barycentric cs:g=0.3,h=0.3,c=0.3) circle (3pt);
\fill[red] (barycentric cs:a=0.25,g=0.25,i=0.25,f=0.25) circle (3pt);
\fill[red] (barycentric cs:h=0.25,g=0.25,i=0.25,j=0.25) circle (3pt);
\fill[red] (barycentric cs:i=0.2,j=0.2,k=0.2,l=0.2,f=0.2) circle
(3pt);
\fill[red] (barycentric cs:f=0.25,l=0.25,m=0.25,e=0.25) circle (3pt);


\fill[red] (4.6,-1) circle (3pt);
\node[right] at (4.7,-1) {$(u_c)_{c\in\mathcal{C}}$};
\path[draw,line width=1mm,vert] (4.5,-2) -- (4.8,-2);
\node[right] at (4.8,-2) {$\Gamma_h(t_k)$};
\end{tikzpicture}
\caption{Domain $\Omega(t_k)$ covered by a polygonal mesh and vector-valued degrees of freedom for the displacement.}
\label{fig:dofs bis}
\end{center}
\end{figure}

Vector-valued volumetric degrees of freedom (dofs) for a generic 
displacement field $(v_c)_{c\in\mathcal{C}}\in\mathbb{R}^{d\#(\mathcal{C})}$ are placed at the barycentre of every mesh cell $c\in\mathcal{C}$. We use the compact notation $v_h:=(v_c)_{c\in\mathcal{C}}$ for the collection of all the cell dofs and we write $v_h \in V_h:=\mathbb{R}^{d\#(\mathcal{C})}$.
Figure \ref{fig:dofs bis} illustrates the position of the displacement
dofs.

\subsection{Discrete bilinear and linear forms}

The discrete stiffness bilinear form hinges on a reconstruction operator that provides a displacement value at every mesh facet by an interpolation formula from neighbouring cell
dofs. 
Specifically, using the cell dofs of $v_h\in V_h$ and the Dirichlet boundary conditions, we reconstruct a collection of displacements $v_{\mathcal{F}}:=(v_F)_{F\in\mathcal{F}(t_k)}\in
\mathbb{R}^{d\#(\mathcal{F}(t_k))}$ on all the mesh facets. 
The reconstruction operator is
denoted $\mathcal{R}(t_k;\cdot)$ and we write
\begin{equation}
v_{\mathcal{F}} := \mathcal{R}(t_k;v_h) \in \mathbb{R}^{d\#(\mathcal{F}(t_k))}.
\end{equation}
The reconstruction operator depends on $t_k$ because of the connectivity modifications due to the crack propagation.

Let us first describe the reconstruction operator on boundary facets.
Let $F \in \mathcal{F}^b_D$ be a Dirichlet boundary facet. Then the reconstruction is simply defined by
evaluating the Dirichlet boundary condition at $\mathbf{x}_F$.
Let $F \in \mathcal{F}^b_N(t_k)$ be a Neumann boundary facet.
The main idea to define $v_F$ is to use a barycentric combination of the cell dofs close to $F$. A similar idea has been considered for finite volume methods in
\cite[Sec. 2.2]{eymard2009discretization} and for cell-centered Galerkin methods in \cite{di2012cell}. We thus select a subset of neighboring cell dofs of $F$, say
$\mathcal{I}_F \subset \mathcal{C}$, and set
\begin{equation}\label{eq:barycentric_crack}
  v_F :=
  \sum_{i\in\mathcal{I}_F}{\alpha_i(\mathbf{x}_F)v_i},
\end{equation}
where the $v_i$'s are the dofs of $v_h$ and 
the coefficients $\alpha_i(\mathbf{x}_F)$ are the barycentric coordinates of the facet barycenter 
$\mathbf{x}_F$ in terms of the selected positions of the dofs. For this construction to be 
meaningful, all the points associated with the selected dofs must not lie on the same line,
so that, in particular, the cardinality of $\mathcal{I}_F$ is at least $(d+1)=3$.

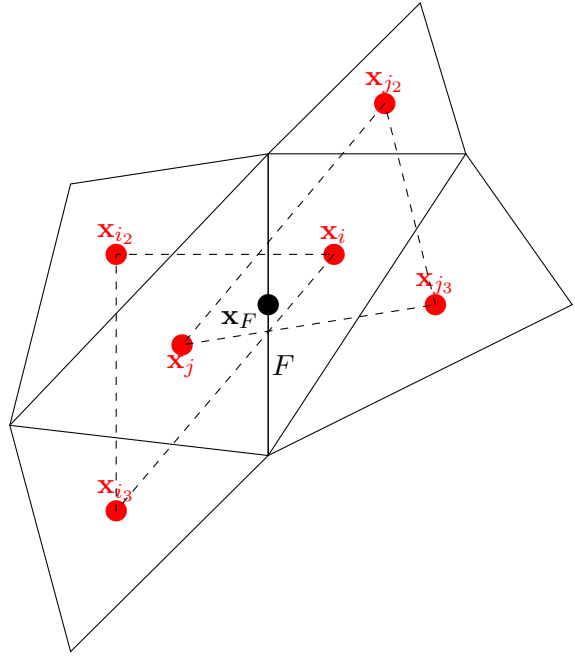
\begin{figure} [!htp]
\begin{center}
\begin{tikzpicture} [scale = 2.]
\coordinate (a) at (0,-1);
\coordinate (b) at (0,1);
\coordinate (c) at (1.3,1);
\coordinate (abc) at (barycentric cs:a=0.3,b=0.3,c=0.3) {};
\path[draw] (a)-- (b)-- (c) -- (a);
\fill[red] (abc) circle (2pt);
\draw (0.43,0.34) node[above,red]{$\mathbf{x}_{i}$};

\coordinate (d) at (-1.7,-0.8);
\coordinate (abd) at (barycentric cs:a=0.3,b=0.3,d=0.3) {};
\path[draw] (a)-- (b)-- (d) -- cycle;
\fill[red] (abd) circle (2pt);
\draw (-0.57,-0.27) node[below, red]{$\mathbf{x}_{j}$};

\coordinate (e) at (-1.3,0.8);
\coordinate (bde) at (barycentric cs:b=0.3,e=0.3,d=0.3) {};
\path[draw] (d)-- (e) -- (b);
\fill[red] (bde) circle (2pt);
\draw (bde) node[above,red]{$\mathbf{x}_{i_2}$};

\coordinate (f) at (-1.3,-2.3);
\coordinate (adf) at (barycentric cs:a=0.3,f=0.3,d=0.3) {};
\path[draw] (d)-- (f) -- (a);
\fill[red] (adf) circle (2pt);
\draw (adf) node[above,red]{$\mathbf{x}_{i_3}$};

\coordinate (g) at (1.,2.);
\coordinate (bcg) at (barycentric cs:b=0.3,c=0.3,g=0.3) {};
\path[draw] (b)-- (g) -- (c);
\fill[red] (bcg) circle (2pt);
\draw (bcg) node[above,red]{$\mathbf{x}_{j_2}$};

\coordinate (h) at (2,0);
\coordinate (ach) at (barycentric cs:a=0.3,c=0.3,h=0.3) {};
\path[draw] (a)-- (h) -- (c);
\fill[red] (ach) circle (2pt);
\draw (ach) node[above,red]{$\mathbf{x}_{j_3}$};

\draw (0.1, -0.4) node {$F$};
\coordinate (ab) at (barycentric cs:a=0.5,b=0.5) {};
\fill (ab) circle (2pt);
\draw (0.,-0.1) node[left] {$\mathbf{x}_F$};
\path[draw, dashed] (abc) -- (adf)-- (bde) -- cycle;
\path[draw, dashed] (ach) -- (bcg)-- (abd) -- cycle;
\end{tikzpicture}
\end{center}
\caption{Dofs associated with the interior facet $F$ used in the reconstruction.}
\label{fig:facet symmetric reconstruction}
\end{figure}

Let us then describe the reconstruction for an inner facet $F\in\mathcal{F}^i(t_k)$.
We use a reconstruction similar to the one presented above except that the two cells sharing the inner facet $F$ play symmetric roles. We refer to this construction as symmetric reconstruction. Specifically, let $c_+$ and $c_-$ be the two cells sharing the inner facet $F\in\mathcal{F}^i(t_k)$. Then, we select $\mathcal{I}_-$ (resp. $\mathcal{I}_+$) as being composed of the cell $c_+$ (resp. $c_-$) and of all the other cells sharing an inner facet with $c_-$ (resp. $c_+$). Notice that these two sets are disjoint. We then set
\begin{equation}\label{eq:barycentric_symm}
  v_F := \frac12  
  \sum_{i\in\mathcal{I}_-\cup\mathcal{I}_+} \alpha_i(\mathbf{x}_F) v_i,
\end{equation}
so that, in the case of a simplicial mesh, $2(d+1)$ dofs are used for the reconstruction (always including $c_-$ and $c_+$). Note that $\sum_{i\in\mathcal{I}_-}\alpha_i(\mathbf{x}_F)=\sum_{i\in \mathcal{I}_+}\alpha_i(\mathbf{x}_F)=1$ here. 
Figure \ref{fig:facet symmetric reconstruction} presents an example where $c_-=c_i$, $c_+=c_j$, $\mathcal{I}_-=\{j,j_2,j_3\}$ and $\mathcal{I}_+=\{i,i_2,i_3\}$.

Having defined the reconstructed facet displacements, it is now possible 
to devise a discrete $\mathbb{R}^{d\times d}$-valued
piecewise-constant gradient field for the displacement
that we write $G_{\mathcal{C}}(v_{\mathcal{F}}) := (G_{c}(v_{\mathcal{F}}))_{c\in\mathcal{C}}
\in \mathbb{R}^{d^2\#(\mathcal{C})}$. Specifically, we set in every
mesh cell $c\in\mathcal{C}$,
\begin{equation}
\label{eq:gradient reconstruction bis}    
G_c(v_{\mathcal{F}}) := \sum_{F \in \partial c} \frac{|F|}{|c|} v_F \otimes n_{F,c},
\end{equation}
where the summation is over the facets $F$ of $c$ and $n_{F,c}$ is the outward 
normal to $c$ on $F$. Note that \eqref{eq:gradient reconstruction bis} is motivated by 
a Stokes formula and that for all $v_h\in V_h$, we have
\begin{equation}
G_c(\mathcal{R}(t_k;v_h)) = \sum_{F \in \partial c} \frac{|F|}{|c|} (\mathcal{R}(t_k;v_h)_F -
v_c) \otimes n_{F,c},
\end{equation} 
since $\sum_{F \in \partial c}|F|n_{F,c}=0$.
We define a constant
linearized strain tensor in every mesh cell $c\in\mathcal{C}$ such that 
\begin{equation}
\varepsilon_{c}(v_{\mathcal{F}}) := \frac{1}{2}(G_{c}(v_{\mathcal{F}})+
G_{c}(v_{\mathcal{F}})\Tr)\in \mathbb{R}^{d\times d},
\end{equation} 
and a constant stress tensor in every mesh cell $c\in\mathcal{C}$ such that 
\begin{equation}
\Sigma_{c}(v_{\mathcal{F}}) := \mathbb{C}:\varepsilon_{c}(v_{\mathcal{F}}) 
\in \mathbb{R}^{d\times d}.
\end{equation}
Finally, we define an additional reconstruction that is used to formulate the stabilization bilinear form in the discrete problem (see below). This operator is a cellwise nonconforming $P^1$ reconstruction $\mathfrak{R}_{c}$ defined for all $c \in \mathcal{C}$ by
\begin{equation}
\label{eq:DG reconstruction}
\mathfrak{R}_{c}(t_k;v_h)(\mathbf{x}) := v_c + G_c(\mathcal{R}(t_k;v_h)) \cdot (\mathbf{x} - \mathbf{x}_c), \qquad \forall \mathbf{x} \in c.
\end{equation}

\subsection{Discrete problem}

We set
\begin{equation} \label{eq:discrete_spaces} \left\{ \begin{aligned}
V_{hD}(t_k) &:=\{v_h\in V_h
\ | \ \mathcal{R}(t_k; v_h)_F = u_D(t_k;\mathbf{x}_F), \ \forall F\subset \partial\Omega_D\}, \quad \forall k=0,\cdots,K,
\\
V_{h0}(t_k) &:=\{v_h\in V_h
\ | \ \mathcal{R}(t_k; v_h)_F = 0, \ \forall F \subset\partial\Omega_D\}, \quad \forall k=0,\cdots,K.
\end{aligned}\right. \end{equation}
The discrete stiffness bilinear form is such that
for all $(v_h,\tilde{v}_h)\in V_{hD}(t_k)\times V_{h0}(t_k)$
(compare with (\ref{eq:bilinear form crack}))
\begin{equation}
a_h(t_k;v_h,\tilde{v}_h) :=
\sum_{c\in\mathcal{C}}|c|\varepsilon_{c}(\mathcal{R}(t_k;v_h)) : \mathbb{C}
: \varepsilon_{c}(\mathcal{R}(t_k;\tilde{v}_h)) + s_h(t_k;v_h,\tilde{v}_h),
\end{equation}
where the stabilization bilinear form $s_h$ is intended to
render $a_h$ coercive and is defined as
\begin{equation}
\label{eq:penalty term crack}
s_h(t_k;v_h,\tilde{v}_h) = \sum_{F \in \mathcal{F}^i(t_k)}  \frac{2\mu}{h_F} |F| [\mathfrak{R}(t_k;v_h)]_F \cdot [ \mathfrak{R}(t_k;\tilde{v}_h)]_F + \sum_{F \in \mathcal{F}^b_D}  \frac{2\mu}{h_F} |F| [\mathfrak{R}(t_k;v_h)]_F \cdot [ \mathfrak{R}(t_k;\tilde{v}_h)]_F,
\end{equation} 
where $h_F$ is the diameter of the facet $F\in\mathcal{F}(t_k)$.
For an interior facet $F\in\mathcal{F}^i(t_k)$, writing $c_-$
and $c_+$ the two mesh cells sharing $F$, i.e., $F=\partial c_-\cap \partial c_+$, 
and orienting $F$ by the unit normal
vector $n_F$ pointing from $c_-$ to $c_+$, the jump of $\mathfrak{R}(t_k;v_h)$ across $F$ is defined as
\begin{equation}
[\mathfrak{R}(t_k;v_h)]_F :=
\mathfrak{R}_{c_-}(t_k;v_h)(\mathbf{x}_F) - \mathfrak{R}_{c_+}(t_k;v_h)(\mathbf{x}_F).
\end{equation}
The sign of the jump is irrelevant in what follows. The role of the 
summation over the interior facets in~\eqref{eq:penalty term crack}
is to penalize the jumps of the cell reconstruction $\mathfrak{R}$ across the interior facets. 
For a Dirichlet boundary facet $F\in\mathcal{F}^b_D$, we denote $c_-$ the unique mesh cell 
containing $F$, we orient $F$ by the unit normal vector $n_F:=n_{c_-}$ which
points outward $\Omega$, and we define
\begin{equation}
[\mathfrak{R}(t_k;v_h)]_F := \mathcal{R}(t_k;v_h)_F - \mathfrak{R}_{c_-}(t_k;v_h)(\mathbf{x}_F).
\end{equation}
Let us recall that for $u_h \in V_{hD}(t_k)$, $\mathcal{R}(t_k;u_h)_F = u_D(t_k;\mathbf{x}_F)$ and for $v_h \in V_{h0}(t_k)$, $\mathcal{R}(t_k;v_h)_F = 0$.
The role of the summation over the Dirichlet boundary facets 
in~\eqref{eq:penalty term crack} is to
penalize the jumps between the cell reconstruction $\mathfrak{R}$ and the value interpolated in the Dirichlet boundary facets.
The bilinear form $s_h$ is classical in the context of discontinuous Galerkin methods (see \cite{arnold1982interior,ern_discontinuous} for instance, see also \cite{di2012cell} for cell-centred Galerkin methods).
It is possible to replace the coefficient $2\mu$ in \eqref{eq:penalty term crack} by $\beta\mu$ with a user-dependent dimensionless parameter $\beta$ of order unity. The numerical experiments reported in \cite{MaraDEM} indicate that this choice has a marginal influence on the results.

\subsection{Verification test case}
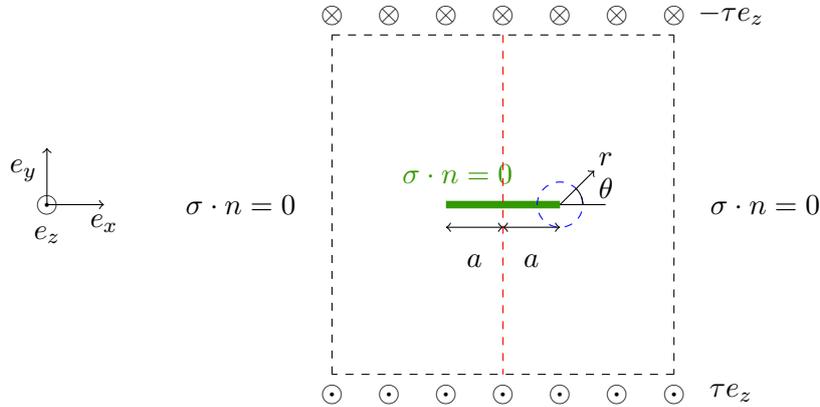
\begin{figure} [!htp]
\begin{center}
\begin{tikzpicture} [scale = 1.5]
\pgfmathsetmacro{\L}{3}
\pgfmathsetmacro{\a}{0.5}
\draw[dashed]  (\L/2,-\L/2) -- (-\L/2,-\L/2);
\draw[dashed]  (-\L/2,\L/2) -- (\L/2,\L/2);
\draw[dashed]  (\L/2,\L/2) -- (\L/2,-\L/2);
\draw[dashed]  (-\L/2,-\L/2) -- (-\L/2,\L/2);

\draw[vert,line width=1mm] (-\a,0) -- (\a,0);
\draw[<->] (0,-0.2) -- (\a,-0.2);
\draw (\a/2,-0.5) node{$a$};
\draw[<->] (0,-0.2) -- (-\a,-0.2);
\draw (-\a/2,-0.5) node{$a$};
\draw[vert,above] (-0.4,0.1) node{$\sigma \cdot n = 0$};
\draw[blue,dashed] (\a+0.2,0.) arc (0:360:0.2);

\draw[-] (\a,0) -- (\a+0.4,0);
\draw (\a+0.2,0) arc (0:45:0.2);
\draw (\a+0.25,0.15) node[right] {$\theta$};
\draw[->] (\a,0) -- (\a+0.3,0.3);
\draw (\a+0.25,0.4) node[right] {$r$};

\draw[dashed,red] (0,\L/2) -- (0,-\L/2);

\draw[above] (0,\L/2) node {$\otimes$};
\draw[above] (1,\L/2) node {$\otimes$};
\draw[above] (0.5,\L/2) node {$\otimes$};
\draw[above] (1.5,\L/2) node {$\otimes$};
\draw[above] (-0.5,\L/2) node {$\otimes$};
\draw[above] (-1.5,\L/2) node {$\otimes$};
\draw[above] (-1,\L/2) node {$\otimes$};
\draw (2.,\L/2) node[above] {$-\tau e_z$};

\draw[below] (0,-\L/2) node {$\odot$};
\draw[below] (1,-\L/2) node {$\odot$};
\draw[below] (0.5,-\L/2) node {$\odot$};
\draw[below] (1.5,-\L/2) node {$\odot$};
\draw[below] (-0.5,-\L/2) node {$\odot$};
\draw[below] (-1.5,-\L/2) node {$\odot$};
\draw[below] (-1,-\L/2) node {$\odot$};
\draw (2.,-\L/2) node[below] {$\tau e_z$};

\draw (-2.3,0) node{$\sigma \cdot n = 0$};
\draw (2.3,0.) node {$\sigma \cdot n = 0$};

\draw [->] (-4,0) -- (-3.5,0);
\draw (-3.5,0) node[below]{$e_x$};
\draw [->] (-4,0) -- (-4,0.5);
\draw (-4,0.3) node[left]{$e_y$};
\draw (-4,0) node {$\odot$};
\draw (-4,-0.1) node[below]{$e_z$};
\end{tikzpicture}
\caption{Sketch of the antiplane shear experiment in an infinite plate.}
\label{fig:crack in antiplane shear}
\end{center}
\end{figure}

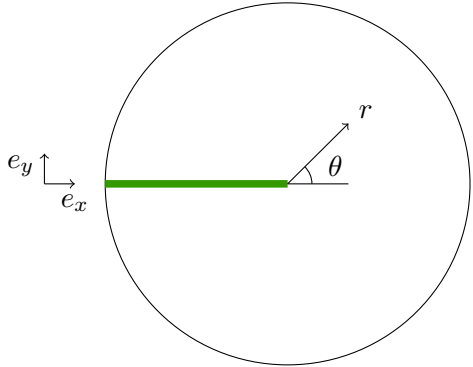
\begin{figure} [!htp]
\begin{center}
\begin{tikzpicture} [scale = 0.8]
\pgfmathsetmacro{\R}{3}
\pgfmathsetmacro{\a}{0}
\draw (\R,0) arc (0:360:\R);
\draw[vert,line width=1mm] (-\R,0) -- (0,0);

\draw[-] (0,0) -- (1,0);
\draw (\a+0.4,0) arc (0:45:0.4);
\draw (\a+0.5,0.3) node[right] {$\theta$};
\draw[->] (0,0) -- (1,1);
\draw (1,1.2) node[right] {$r$};

\draw [->] (-4,0) -- (-3.5,0);
\draw (-3.5,0) node[below]{$e_x$};
\draw [->] (-4,0) -- (-4,0.5);
\draw (-4,0.3) node[left]{$e_y$};
\end{tikzpicture}
\caption{Sketch of the meshed ball around the crack tip.}
\label{fig:ball around crack}
\end{center}
\end{figure}

This section presents a verification test case related to the convergence rate with a singularity at the crack tip.
The crack does not propagate, i.e., we consider a steady setting using the above discrete stiffness bilinear form and load linear form.
The convergence rate of the method in the presence of a singularity is tested in the case of an infinite plate under mode 3 loading at infinity as presented in Figure \ref{fig:crack in antiplane shear}.
A convergence rate of $O(h^{\frac12})$, similar to that obtained with Lagrange $P^1$ finite elements, is expected.
The reference solution, close to the crack tip ($\frac{r}{a} \ll 1$), reads in polar coordinates \cite[p. 28]{kuna2013finite}:
\begin{equation}
\label{eq:ref solution antiplane}
u(r,\theta) = \frac{2 \tau}{\mu} \sqrt{\frac{ar}{2}} \sin \left(\frac{\theta}{2} \right) e_z,
\end{equation}
where $\tau$ is the modulus of the antiplane shear stress imposed at infinity. The displacement defined in (\ref{eq:ref solution antiplane}) verifies the statics equation in a strong form since $\mathrm{div}(u) = 0$. The stresses are
\begin{equation}
\label{eq: stresses antiplane}
\sigma(r,\theta) = \tau \sqrt{\frac{a}{2 r}} \left[\sin \left( \frac{\theta}{2} \right) e_r - \cos \left( \frac{\theta}{2} \right) e_\theta \right] \otimes e_z.
\end{equation}
The domain shown in Figure \ref{fig:crack in antiplane shear} being symmetric with respect to the red dashed line, only its right part is considered. As the analytical solution (\ref{eq:ref solution antiplane}) is only valid close to the crack tip, a small ball around the crack tip, which corresponds to the green dashed circle in Figure \ref{fig:crack in antiplane shear}, is meshed. The setting is presented in Figure \ref{fig:ball around crack}. The convergence towards the analytical solution is checked on the meshed ball with the reference solution imposed as Dirichlet boundary condition over the whole boundary including the crack lips.
The results of the computation, which are reported in Table \ref{tab:convergence antiplane}, corroborate an $O(h^{\frac12})$ convergence rate in the energy-norm, as expected.
\begin{table}[!htp]
\begin{center}
   \begin{tabular}{ | c | c | c | c | c |}
     \hline
     nb dofs & $ \Vert u - \mathfrak{R}(u_h) \Vert_{L^2}$ & Convergence rate & $\Vert \nabla u - G_h(u_h) \Vert_{L^2}$ & Convergence rate  \\ \hline
     $496$ &  5.84e-05 & - & 1.22e-01 & - \\ \hline
     $1,880$ & 1.77e-05 & $1.80$  & 8.16e-02 & $0.57$ \\ \hline
     $7,312$ &  5.76e-06 & $1.65$  & 5.66e-02 & $0.50$  \\ \hline
     $28,832$ & 1.96e-06 & $1.57$ & 3.95e-02 & $0.50$ \\ \hline
     $114,496$ & 6.83e-07 & $1.53$ & 2.78e-02 & $0.50$ \\ \hline
   \end{tabular}
   \caption{Number of dofs, $L^2$-error and convergence rate, $L^2$-error on the gradient and convergence rate.}
   \label{tab:convergence antiplane}
\end{center}
\end{table}
We also observe an $O(h^{\frac32})$ convergence rate in the $L^2$-norm.
The convergence rates are evaluated as
\begin{equation}
\text{order} = d \log\left(\frac{e_1}{e_2}\right)  \left(\log\left( \frac{n_2}{n_1} \right) \right)^{-1},
\end{equation}
where $e_1,e_2$ denote the errors on the computations with mesh sizes $h_1,h_2$ and the number of dofs $n_1,n_2$.

\section{Quasi-static crack propagation}
\label{sec:quasi-static cracking}
In this section, we formulate the discrete problem for quasi-static crack propagation. The space discretization is achieved by means of the variational DEM scheme presented in the previous section.
At every pseudo-time node $t_k$, the problem is solved iteratively with inner iterations enumerated by $m \in \{0, \ldots, M\}$. Since the crack can change at each inner iteration, we use the notation $\Gamma_h(t_{k,m})$ for the crack and the notation $\mathcal{F}^i(t_{k,m})$ and $\mathcal{F}^b(t_{k,m})$ for the partition of the mesh facets at the inner iteration $m$, with the facets located in the crack collected in the subset $\mathcal{F}^\Gamma(t_{k,m})$. 

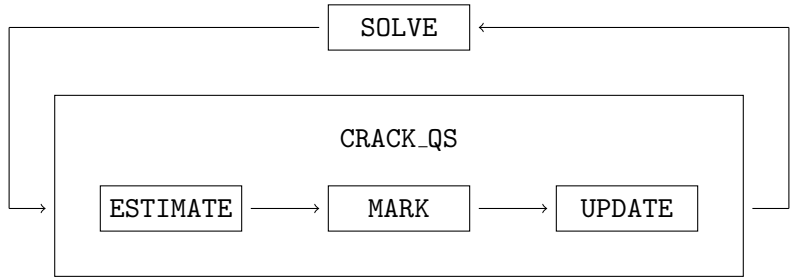
\begin{figure} [!htp]
\begin{center}
\begin{tikzpicture} [scale = 0.6]
\pgfmathsetmacro{\L}{3.1}
\pgfmathsetmacro{\H}{1}
\pgfmathsetmacro{\po}{5}
\pgfmathsetmacro{\pt}{10}
\pgfmathsetmacro{\dl}{0.2}
\pgfmathsetmacro{\bigH}{3}

\draw (\L+\po, 5) -- (\L+\po, 4) -- (\po, 4) -- (\po, 5) -- cycle;
\draw (\po+\L/2, 4.5) node {\texttt{SOLVE}};
\draw[<-] (\po+\L+\dl, 4.5) -- (\L+\pt+2, 4.5);
\draw (\L+\pt+2, \H/2) -- (\L+\pt+2, 4.5);
\draw (\L+\pt+1+\dl, \H/2) -- (\L+\pt+2, \H/2);
\draw (-2, 4.5) -- (\po-\dl, 4.5);
\draw (-2, \H/2) -- (-2, 4.5);
\draw[->] (-2, \H/2) -- (-1-\dl, \H/2);

\draw (\L, \H) -- (\L, 0) -- (0, 0) -- (0, \H) -- cycle;
\draw[->] (\L+\dl, \H/2) -- (\po-\dl, \H/2);
\draw (\L/2, \H/2) node {\texttt{ESTIMATE}};

\draw (\L+\po, \H) -- (\L+\po, 0) -- (\po, 0) -- (\po, \H) -- cycle;
\draw[->] (\po+\L+\dl, \H/2) -- (\pt-\dl, \H/2);
\draw (\po+\L/2, \H/2) node {\texttt{MARK}};

\draw (\L+\pt, \H) -- (\L+\pt, 0) -- (\pt, 0) -- (\pt, \H) -- cycle;
\draw (\pt+\L/2, \H/2) node {\texttt{UPDATE}};

\draw (\L+\pt+1, \bigH) -- (\L+\pt+1, -1) -- (-1, -1) -- (-1, \bigH) -- cycle;
\draw (\L/2+\pt/2, \bigH-1) node {\texttt{CRACK\_QS}};
\end{tikzpicture}
\caption{Details of the procedure \texttt{CRACK\_QS}.}
\label{fig:procedure}
\end{center}
\end{figure}

Each inner iteration consists in two steps. First, freezing the position of the crack, we find the discrete displacement
$u_h(t_{k,m})\in V_{hD}(t_k)$ solving the quasi-static problem $a_h(t_{k,m};u_h(t_{k,m}),\tilde{v}_h) = l_h(t_{k};\tilde{v}_h)$
for all $\tilde{v}_h\in V_{h0}(t_{k})$ (the bilinear form $a_h$ depends on $t_{k,m}$ since the reconstruction operator changes as the crack propagates). Then we use the newly
computed displacement field $u_h(t_{k,m})$ to determine whether crack propagation occurs
and update accordingly the subsets $\mathcal{F}^i(t_{k,m+1})$, 
$\mathcal{F}^b(t_{k,m+1})$, and $\mathcal{F}^\Gamma(t_{k,m+1})$. We iterate this procedure until there is no more crack propagation in the second step.
The inner iteration in the discrete quasi-static crack propagation scheme can thus be summarized as follows:
For all $m \in \{0, \ldots, M \}$,
\begin{equation}
\label{eq:quasi-static discretization}
\left\{
\begin{alignedat}{2}
&\textup{(i)}&\ &u_h(t_{k,m})\in V_{hD}(t_{k}) \ \text{s.t.} \  a_h(t_{k,m};u_h(t_{k,m}),\tilde{v}_h) = l_h(t_{k};\tilde{v}_h),\ \forall \tilde{v}_h\in V_{h0}(t_{k}), \\
&\textup{(ii)}&\ &(\mathcal{F}^{\Gamma}(t_{k,m+1}),\mathcal{F}^b(t_{k,m+1}),\mathcal{F}^i(t_{k,m+1})) = \texttt{CRACK\_QS}(\mathcal{F}^{\Gamma}(t_{k,m}),\mathcal{F}^b(t_{k,m}),\mathcal{F}^i(t_{k,m}),u_h(t_{k,m})).
\end{alignedat}
\right.
\end{equation}
The rest of this section is devoted to the description of the procedure \texttt{CRACK\_QS}.
This procedure consists in the three consecutive steps outlined in Figure \ref{fig:procedure}. 
The first step involves the procedure \texttt{ESTIMATE} which considers all the vertices of $\mathcal{F}^{\Gamma}(t_{k,m})$ and computes for each of these vertices an approximate energy release rate. The second step involves the procedure \texttt{MARK} which flags among all the inner facets sharing a vertex with an energy release rate larger than the maximum value $\mathcal{G}_c$ the facet that will indeed break.
The selection is made by using a discrete kinking criterion.
The last step uses the procedure \texttt{UPDATE} and simply consists in updating the data structure according to the crack propagation. The procedure is repeated from the recomputation of the solution of the first line of Equation \eqref{eq:quasi-static discretization} until no facet is marked in the procedure \texttt{MARK}.

\subsection{Procedure \texttt{ESTIMATE}}
\label{sec:estimation of SIF}
Let $\mathcal{V}^{\Gamma}(t_{k,m})$ be the set of all vertices in the crack $\Gamma(t_{k,m})$.
The procedure \texttt{ESTIMATE} computes an approximate energy release rate $\mathcal{G}_h(v)$ for all $v \in \mathcal{V}^{\Gamma}(t_{k,m})$.
Let $\mathcal{F}^{\Gamma}_\mathbf{v}(t_{k,m})$ be the set of cracked facets sharing a vertex $\mathbf{v} \in \mathcal{V}^{\Gamma}(t_{k,m})$. (The set $\mathcal{F}^{\Gamma}_\mathbf{v}(t_{k,m})$ reduces to a single facet if $\mathbf{v}$ is the crack tip.)
Let $\mathcal{F}^i_\mathbf{v}(t_{k,m})$ be the set of inner facets sharing a vertex $\mathbf{v} \in \mathcal{V}^{\Gamma}(t_{k,m})$.
An approximate energy release rate for
the vertex $\mathbf{v} \in \mathcal{V}^{\Gamma}(t_{k,m})$ is evaluated as
\begin{equation}
\mathcal{G}_h(\mathbf{v}) := \max_{F \in \mathcal{F}^{\Gamma}_\mathbf{v}(t_{k,m})} \max_{F' \in \mathcal{F}^i_\mathbf{v}(t_{k,m})} \pi n_F \cdot \{\Sigma_h(t_{k,m})\}_{F} \cdot [u_h(t_{k,m})]_{F'},
\end{equation}
where $[u_h]_F := u_{c_-} - u_{c_+}$, $\{\Sigma_h\}_F := \frac12 (\Sigma_{c_-} + \Sigma_{c_+})$, and $n_{F}$ is the normal vector to $F$ pointing from $c_-$ to $c_+$.
This expression is rooted in the fact that the elastic energy contained in a facet $F$ writes $\frac12 n_F \cdot \{\Sigma_h(t_{k,m})\}_{F} \cdot [u_h(t_{k,m})]_{F} |F| $ as motivated in \cite{MaraDEM}. The factor $\pi$ comes from the fact that the density of elastic energy per facet must be multiplied by $2\pi$ to take into account the surface created by cracking (see \cite[p. 48]{kuna2013finite}). This is linked to the concept of the crack closure integral.
The output of the procedure \texttt{ESTIMATE} is the collection of approximate energy release rates $\{\mathcal{G}_h(\mathbf{v})\}_{\mathbf{v} \in \mathcal{V}^{\Gamma}(t_{k,m})}$.

\subsection{Procedure \texttt{MARK}}
\label{sec:kinking criterion}
The goal of the procedure \texttt{MARK} is to identify the unique inner facet $\mathfrak{F} \in \mathcal{F}^i(t_{k,m})$ through which the crack will propagate.
The criterion is based on an adaptation of the maximisation of the strain energy density which was introduced in \cite{sih1974strain}.
The vertices of $\mathcal{V}^{\Gamma}(t_{k,m})$ are ordered as they break during a computation and we select the last $N$ vertices in $\mathcal{V}^{\Gamma}(t_{k,m})$ to define the subset $\mathcal{V}^{\Gamma}_N(t_{k,m})$. The integer parameter $N$ is set to $N=6$ in our computations; this choice gives satisfactory results while avoiding excessive branching of the crack path. 
Finally, we select the vertices in $\mathcal{V}^{\Gamma}_N(t_{k,m})$ whose approximate energy release rate is larger than the material parameter $\mathcal{G}_c$:
\begin{equation}
\label{eq:set of boundary entities of the crack}
\mathcal{V}^{\Gamma*}_N(t_{k,m}) := \{ \mathbf{v} \in \mathcal{V}^{\Gamma}_N(t_{k,m}), \mathcal{G}_h(\mathbf{v}) \ge \mathcal{G}_c \}.
\end{equation}
Among all $\mathbf{v} \in \mathcal{V}^{\Gamma*}_N(t_{k,m})$, we select the single vertex through which the crack will propagate at $t_{k,m}$ as
\begin{equation}
\mathbf{z} := \mathop{\mathrm{Argmax}} \limits_{\mathbf{v} \in \mathcal{V}^{\Gamma*}_N(t_{k,m})} \mathcal{G}_h(\mathbf{v}).
\end{equation}
If there is more than one maximizer, one is picked randomly. Note that in most situations, the vertex $\mathbf{z}$ is located at the crack tip.

\begin{figure} [!htp]
\centering
\begin{tikzpicture}[scale=0.8]
\coordinate (b) at (1,-1);
\coordinate (a) at (0,1);
\coordinate (c) at (1,-3);
\coordinate (d) at (0,-5);
\coordinate (e) at (-1,0);
\coordinate (f) at (-2,-3);
\coordinate (g) at (2.7,1.5);
\coordinate (h) at (2.5,-1);
\coordinate (i) at (2.2,-4);
\coordinate (j) at (1.7,-4.6);
\coordinate (k) at (1.8,-1.2);
\coordinate (l) at (0.9,1.5);
\coordinate (m) at (1.7,-4.2);
\coordinate (n) at (0.85,1.15);
\coordinate (l) at (-5,-3);

\path[draw] (c) -- (d);
\draw[red, line width=1mm, opacity=0.5] (a) -- (b) -- (c);
\draw (a) -- (e);
\draw[red, line width=1mm, opacity=0.5] (e) -- (b);
\draw (f) -- (c);
\draw (d) -- (f);
\draw (e) -- (f);
\draw (a) -- (g);
\draw[red, line width=1mm, opacity=0.5] (g) -- (b);
\draw (h) -- (c);
\draw (c) -- (i) -- (d);
\draw (g) -- (h);
\draw[blue, line width = 1mm, opacity=0.5] (h) -- (b);

\draw[vert, line width=1mm, opacity=0.7] (b) -- (f);
\draw[vert, line width=1mm, dashed, opacity=0.7] (f) -- (l);

\fill (b) circle (3pt);
\draw (b) node[above] {$\mathbf{z}$};
\fill (-2.1,-2.9) -- (-1.9,-2.9) -- (-1.9,-3.1) -- (-2.1,-3.1);

\draw[vert, line width = 1mm, opacity=0.7] (5,-1.25) -- (5.5,-1.25);
\draw (6.5,-1.25) node {$\mathcal{F}^{\Gamma}(t_{k,m})$};
\draw[red, line width=1mm, opacity=0.5] (5,-2.75) -- (5.5,-2.75);
\draw (7.5,-2.75) node {$F\in \mathcal{F}^i_\mathbf{z}(t_{k,m}) \setminus \{\mathfrak{F}\}$};
\fill (5.25,0) circle (3pt);
\draw (6.5,0) node {$\mathcal{V}^{\Gamma*}(t_{k,m})$};
\fill (5.2,1) -- (5.4,1) -- (5.4,1.2) -- (5.2,1.2);
\draw (6.5,1.) node {$\mathcal{V}^{\Gamma}(t_{k,m})$};
\draw[blue, line width=1mm, opacity=0.5] (5,-4.25) -- (5.5,-4.25);
\draw (6.,-4.25) node {$\mathfrak{F}$};
\end{tikzpicture}
\caption{Sketch of the discrete sets considered in the procedure \texttt{MARK}.}
\label{fig:discrete sets}
\end{figure}
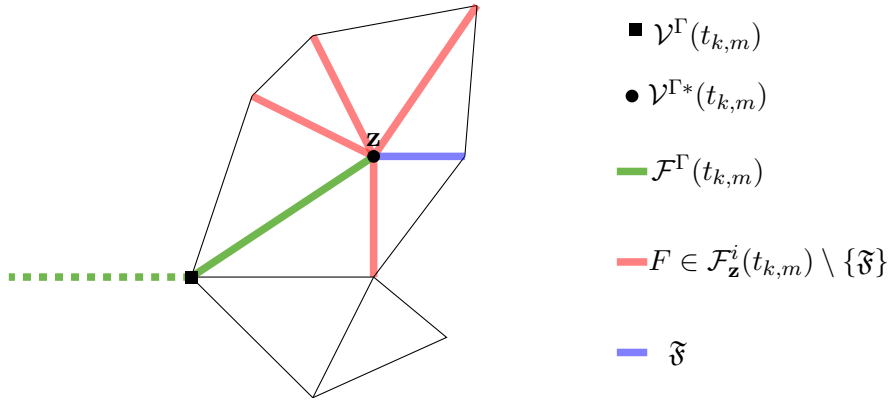

Having selected the vertex $\mathbf{z}$, we now mark one facet $\mathfrak{F} \in \mathcal{F}^i_\mathbf{z}(t_{k,m})$ for cracking.
We impose only one restriction on the selection process of the facet to be broken: we limit the number of facets broken per cell to one.
This limit is justified by the fact that when a facet breaks, the resulting geometric singularity creates very high stresses that lead to breaking the other facets of the cells containing the facet thus creating many fragments. The limitation we impose is to avoid this situation.
The setting is illustrated in Figure \ref{fig:discrete sets}.
The output of the procedure \texttt{MARK} is the facet $\mathfrak{F}$, through which the crack will propagate, defined as
\begin{equation}
\label{eq:kinking criterion}
\mathfrak{F} := \mathop{\mathrm{Argmax}} \limits_{\substack{F \in \mathcal{F}^i_\mathbf{z}(t_{k,m}) \setminus \mathcal{F}^i_\mathcal{C}(t_{k,m})}} 
\frac12 \{\Sigma_h(t_{k,m})\}_F \cdot \{\varepsilon_h(t_{k,m})\}_F,
\end{equation}
where $\mathcal{F}^i_\mathcal{C}(t_{k,m})$ denotes the set of inner facets contained in a cell with one facet already broken.

\subsection{Procedure \texttt{UPDATE}}
\label{sec:rigidity matrix modifications}
The subsets $\mathcal{F}^{\Gamma}(t_{k,m+1})$, $\mathcal{F}^{i}(t_{k,m+1})$, and
$\mathcal{F}^{b}(t_{k,m+1})$ can now be updated as follows:
\begin{equation}
\label{eq:update in the sets of facets}
\left\{
\begin{aligned}
& \mathcal{F}^{\Gamma}(t_{k,m+1})  := \mathcal{F}^{\Gamma}(t_{k,m}) \cup \{\mathfrak{F}\}, \\
& \mathcal{F}^i(t_{k,m+1}) := \mathcal{F}^i(t_{k,m}) \setminus \{\mathfrak{F}\}, \\
& \mathcal{F}^b(t_{k,m+1}) := \mathcal{F}^b(t_{k,m}) \cup  \{\mathfrak{F}_-,\mathfrak{F}_+\},
\end{aligned}
\right.
\end{equation}
where we recall that $\mathfrak{F}_-$ and $\mathfrak{F}_+$ are the same geometric object as the inner facet $\mathfrak{F}$, but are now each one on the boundary of a single mesh cell.

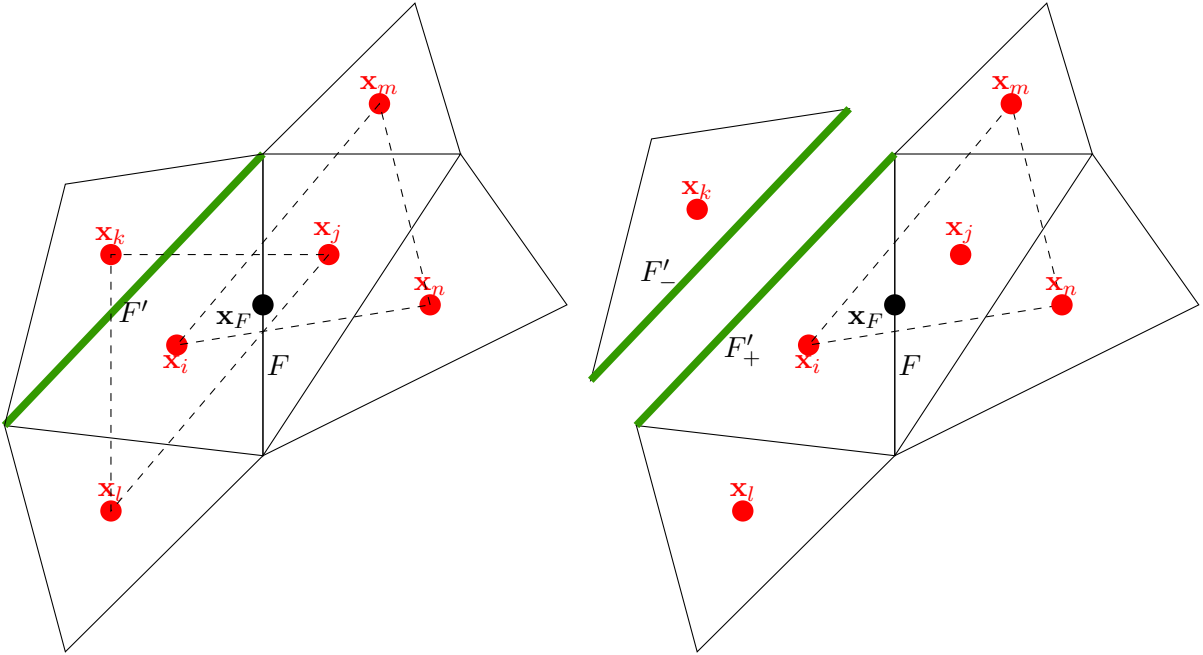
\begin{figure} [!htp]
\centering
\subfloat{
\begin{tikzpicture} [scale = 2.]
\coordinate (a) at (0,-1);
\coordinate (b) at (0,1);
\coordinate (c) at (1.3,1);
\coordinate (abc) at (barycentric cs:a=0.3,b=0.3,c=0.3) {};
\path[draw] (a)-- (b)-- (c) -- cycle;
\fill[red] (abc) circle (2pt);
\draw (0.43,0.34) node[above,red]{$\mathbf{x}_j$};

\coordinate (d) at (-1.7,-0.8);
\coordinate (abd) at (barycentric cs:a=0.3,b=0.3,d=0.3) {};
\path[draw] (a)-- (b)-- (d) -- cycle; 
\fill[red] (abd) circle (2pt);
\draw (-0.57,-0.27) node[below, red]{$\mathbf{x}_i$};
\draw [vert, line width=1mm] (b) -- (d);
\draw (barycentric cs:b=0.5,d=0.5) node[below] {$F'$};

\coordinate (e) at (-1.3,0.8);
\coordinate (bde) at (barycentric cs:b=0.3,e=0.3,d=0.3) {};
\path[draw] (d)-- (e) -- (b);
\fill[red] (bde) circle (2pt);
\draw (bde) node[above,red]{$\mathbf{x}_k$};

\coordinate (f) at (-1.3,-2.3);
\coordinate (adf) at (barycentric cs:a=0.3,f=0.3,d=0.3) {};
\path[draw] (d)-- (f) -- (a);
\fill[red] (adf) circle (2pt);
\draw (adf) node[above,red]{$\mathbf{x}_l$};

\coordinate (g) at (1.,2.);
\coordinate (bcg) at (barycentric cs:b=0.3,c=0.3,g=0.3) {};
\path[draw] (b)-- (g) -- (c);
\fill[red] (bcg) circle (2pt);
\draw (bcg) node[above,red]{$\mathbf{x}_m$};

\coordinate (h) at (2,0);
\coordinate (ach) at (barycentric cs:a=0.3,c=0.3,h=0.3) {};
\path[draw] (a)-- (h) -- (c);
\fill[red] (ach) circle (2pt);
\draw (ach) node[above,red]{$\mathbf{x}_n$};

\draw (0.1, -0.4) node {$F$};
\coordinate (ab) at (barycentric cs:a=0.5,b=0.5) {};
\fill (ab) circle (2pt);
\draw (0.,-0.1) node[left] {$\mathbf{x}_F$};
\path[draw, dashed] (abc) -- (adf)-- (bde) -- cycle;
\path[draw, dashed] (ach) -- (bcg)-- (abd) -- cycle;
\end{tikzpicture}
}
\subfloat{
\begin{tikzpicture} [scale = 2.]
\pgfmathsetmacro{\dl}{0.3}

\coordinate (a) at (0,-1);
\coordinate (b) at (0,1);
\coordinate (c) at (1.3,1);
\coordinate (abc) at (barycentric cs:a=0.3,b=0.3,c=0.3) {};
\path[draw] (a)-- (b)-- (c) -- (a);
\fill[red] (abc) circle (2pt);
\draw (0.43,0.34) node[above,red]{$\mathbf{x}_j$};

\coordinate (d) at (-1.7,-0.8);
\coordinate (abd) at (barycentric cs:a=0.3,b=0.3,d=0.3) {};
\path[draw] (a)-- (b)-- (d) -- cycle;
\fill[red] (abd) circle (2pt);
\draw (-0.57,-0.27) node[below, red]{$\mathbf{x}_i$};
\draw [vert, line width=1mm] (b) -- (d);

\coordinate (ep) at (-1.3-\dl,0.8+\dl);
\coordinate (dp) at (-1.7-\dl,-0.8+\dl);
\coordinate (bp) at (0-\dl,1+\dl);
\coordinate (bdep) at (barycentric cs:bp=0.3,ep=0.3,dp=0.3) {};
\path[draw] (dp)-- (ep) -- (bp) -- cycle;
\fill[red] (bdep) circle (2pt);
\draw (bdep) node[above,red]{$\mathbf{x}_k$};
\draw [vert, line width=1mm] (bp) -- (dp);

\draw (0.1, -0.4) node {$F$};
\coordinate (ab) at (barycentric cs:a=0.5,b=0.5) {};
\fill (ab) circle (2pt);
\draw (0.,-0.1) node[left] {$\mathbf{x}_F$};
\path[draw, dashed] (ach) -- (bcg)-- (abd) -- cycle;

\draw (-1.55, 0.2) node {$F'_-$};
\draw (-1, -0.3) node {$F'_+$};

\coordinate (f) at (-1.3,-2.3);
\coordinate (adf) at (barycentric cs:a=0.3,f=0.3,d=0.3) {};
\path[draw] (d)-- (f) -- (a);
\fill[red] (adf) circle (2pt);
\draw (adf) node[above,red]{$\mathbf{x}_l$};

\coordinate (g) at (1.,2.);
\coordinate (bcg) at (barycentric cs:b=0.3,c=0.3,g=0.3) {};
\path[draw] (b)-- (g) -- (c);
\fill[red] (bcg) circle (2pt);
\draw (bcg) node[above,red]{$\mathbf{x}_m$};

\coordinate (h) at (2,0);
\coordinate (ach) at (barycentric cs:a=0.3,c=0.3,h=0.3) {};
\path[draw] (a)-- (h) -- (c);
\fill[red] (ach) circle (2pt);
\draw (ach) node[above,red]{$\mathbf{x}_n$};

\end{tikzpicture}
}
\caption{Recomputation of the reconstruction stencil associated with the inner facet $F$ after the breaking of the neighbouring inner facet $F'$. Left: reconstruction before cracking. Right: reconstruction after cracking. (The two cells separated by the crack are drawn slightly apart.)}
\label{fig:recomputation of reconstructions}
\end{figure}

\begin{remarque}[Update of $a_h$]
The updates in~\eqref{eq:update in the sets of facets} affect the reconstruction operator used to evaluate the discrete stiffness bilinear form. Figure \ref{fig:recomputation of reconstructions} presents a sketch of an inner facet whose reconstruction has to be recomputed after a neighbouring inner facet breaks. The purpose of recomputing the reconstruction on certain inner facets is to avoid using dof values on both sides of the crack in the same reconstruction.
\end{remarque}

\section{Numerical experiments}
\label{sec:numerics}
Several numerical experiments are presented to show the versatility of the proposed numerical method. The python scripts\protect\footnotemark \ for these numerical experiments use the finite element library FEniCS \cite{LoggMardalEtAl2012a}
\footnotetext{\url{https://github.com/marazzaf/DEM_cracking.git}}
and scipy\protect\footnotemark. Although the proposed method is able to handle polyhedral meshes, our computations only use triangular meshes. This is a consequence of the current restriction of FEniCS to simplicial meshes.

\footnotetext{\url{https://scipy.org/}}

\subsection{Crack speed with prescribed crack path}
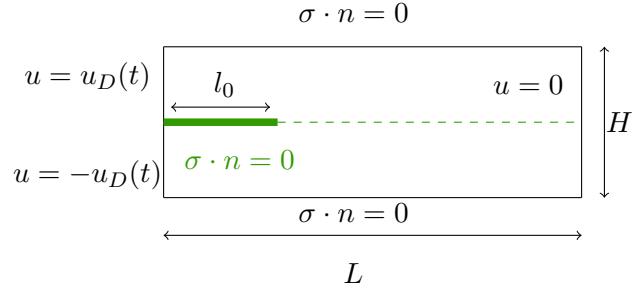
\begin{figure} [!htp]
\centering
\begin{tikzpicture}
\draw (-2.5,1) -- (3,1);
\draw (-2.5,-1) -- (3,-1);
\draw (-2.5,0.) -- (-2.5,1);
\draw (-2.5,0.) -- (-2.5,-1);
\draw (3,1) -- (3,-1);
\draw [<->] (-2.5,-1.5) -- (3,-1.5);
\draw (0,-2) node{$L$};
\draw[vert,line width=1mm] (-2.5,0) -- (-1,0);
\draw (-1.5,-0.5) node[vert] {$\sigma \cdot n = 0$};
\draw [<->] (3.3,-1) -- (3.3,1);
\draw (3.5,0) node{$H$};
\draw[dashed,vert] (-1,0) -- (3,0);
\draw (-1.7,0.5) node{$l_0$};
\draw[<->] (-2.4,0.2) -- (-1.1,0.2);
\draw (0,1.7) node[below]{$\sigma \cdot n = 0$};
\draw (0,-1.2) node{$\sigma \cdot n = 0$};
\draw (2.3,0.5) node {$u = 0$};
\draw (-3.5,-0.7) node {$u = -u_D(t)$};
\draw (-3.5,0.6) node {$u = u_D(t)$};
\end{tikzpicture}
\caption{Crack speed: problem setup.}
\label{fig:sketch antiplane shear constrained}
\end{figure}

\begin{figure}[!htp]
\centering
\subfloat{
\includegraphics[width=0.5\textwidth]{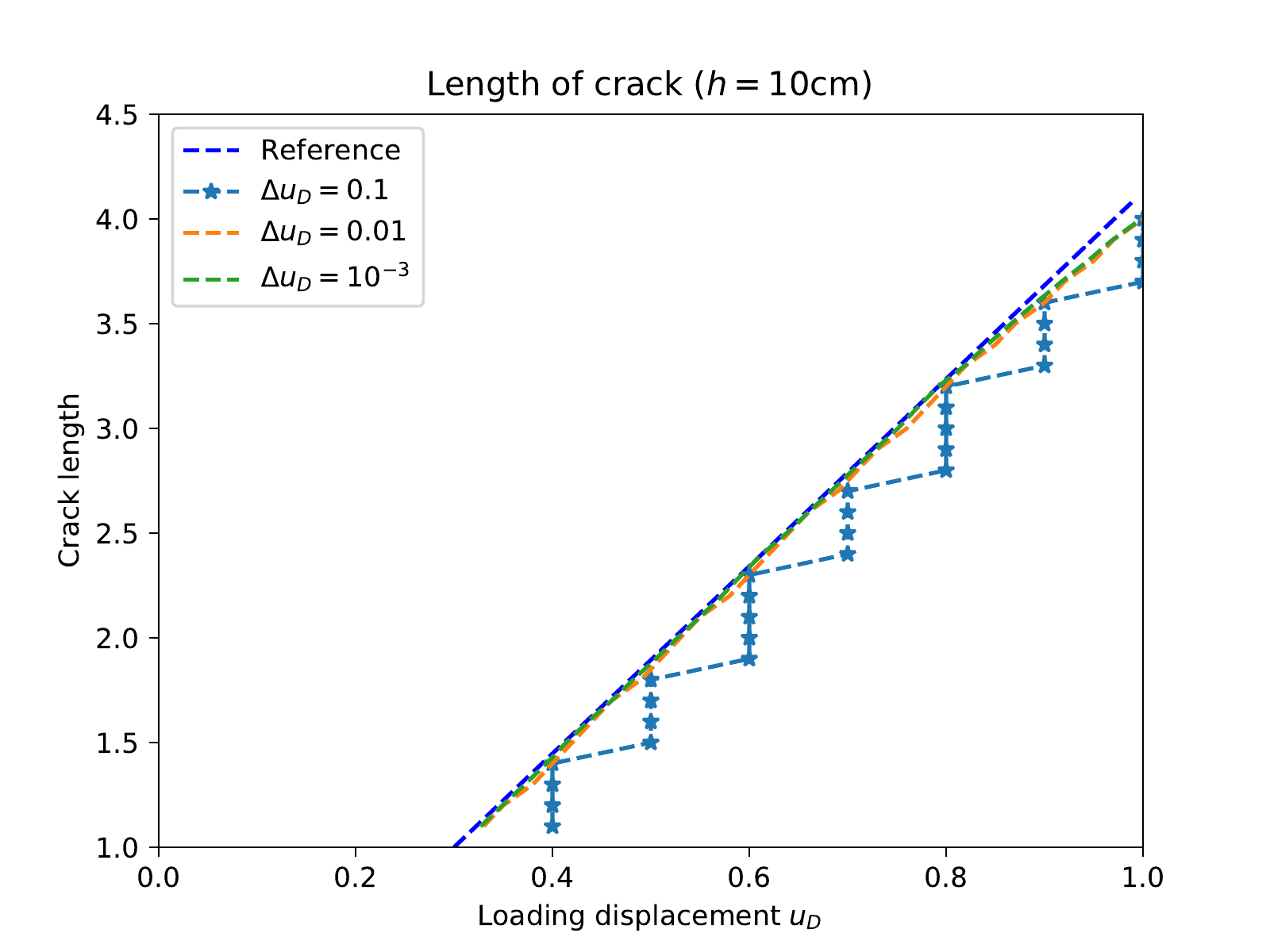}
}
\subfloat{
\includegraphics[width=0.5\textwidth]{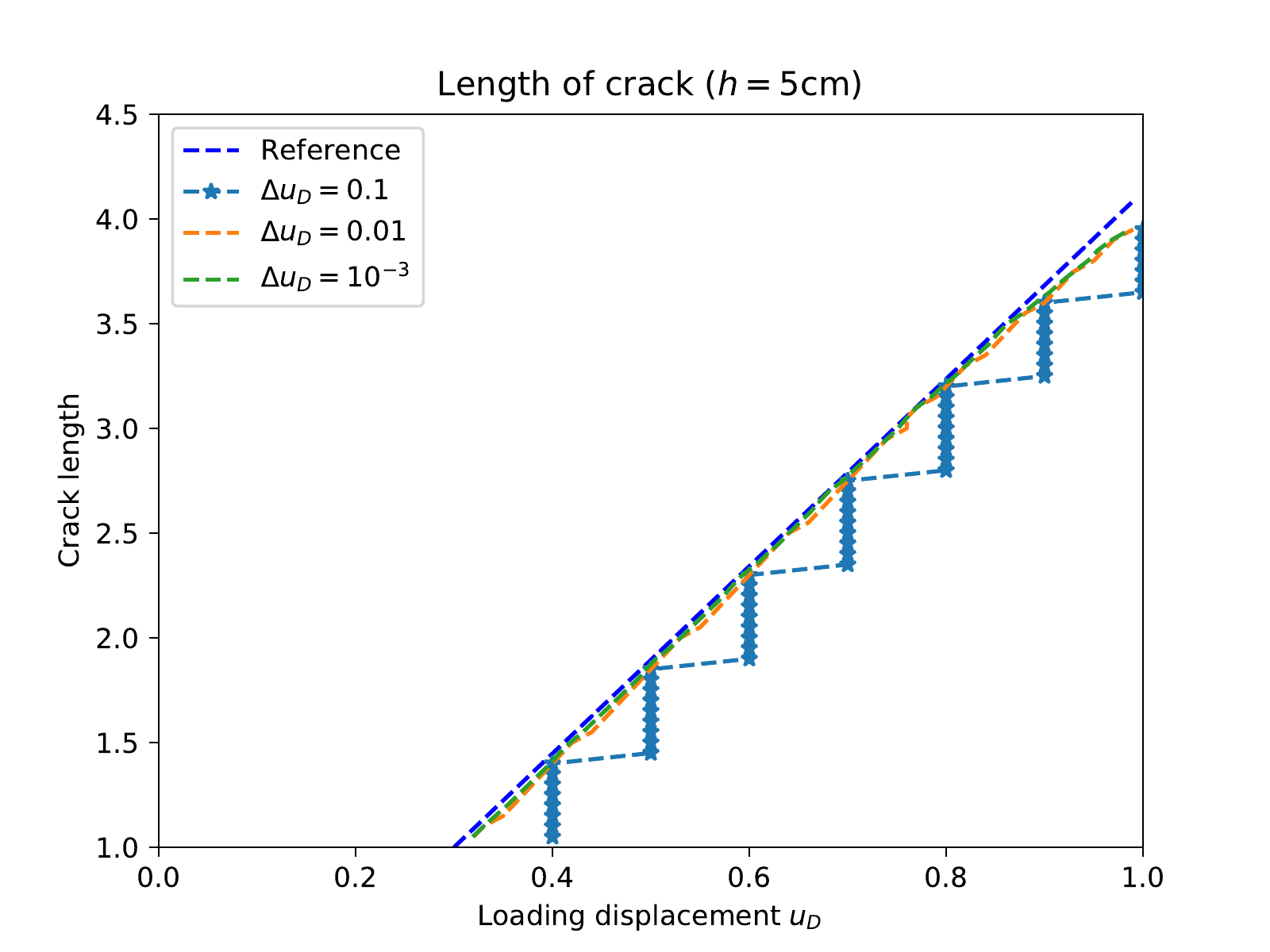}
}
\caption{Crack speed: crack length as a function of loading displacement $u_D$. Left: coarse mesh. Right: fine mesh.}
\label{fig:comparaison antiplane}
\end{figure}

We consider a test case taken from \cite{li2016numerical}. The test case consists of an already cracked plate under antiplane shear loading. 
The crack is forced to propagate along a straight line represented by the dashed line in Figure \ref{fig:sketch antiplane shear constrained}. The goal of this test case is to study the crack propagation velocity.
The dimensions of the plate are $L=5\text{m}$ and $H=1\text{m}$ and the initial length of the crack is $l_0 = 1\text{m}$.
The constant increment in boundary loading is written $\Delta u_D$. The material parameters are $\mu = 0.2 \text{Pa}$ and $\mathcal{G}_c = 0.01 \text{kN/mm}$.
We are interested in the length of the crack with respect to the cumulated boundary loading displacement $u_D$, where the final displacement load is $u_D=1\mathrm{m}$.
The reference solution for the crack speed $S$ with respect to the loading speed, taken from \cite{li2016numerical}, is $\sqrt{\frac{\mu H}{\mathcal{G}_c}} \approx 4.47$. As this solution is only valid when $L \to \infty$, we checked that doubling the length $L$ of the strip did not lead to any significant change in the crack speeds. 
The computations are performed with two structured 2d meshes of triangles with characteristic sizes $h=10$cm and $h=5$cm. Various values of $\Delta u_D$ are used in the two computations.
Figure \ref{fig:comparaison antiplane} reports the crack length as a function of the cumulated loading displacement $u_D$.
One can see that the results with the two meshes are very similar. 
The results with $\Delta u_D = 10^{-3}$m and $\Delta u_D = 10^{-2}$m are very similar and are in agreement with the analytical solution. For these two values, $\frac{\Delta u_D}{h}$ is less than $0.5$, so that the increment in the imposed Dirichlet condition is smaller than the mesh size. This is not the case for $\Delta u_D=0.1$m.
The different aspect of the curves for $\Delta u_D=0.1$m is explained by the fact that as $\Delta u_D$ is large in that case, a large number of facets can break at some of the displacement increments, thus leading to this staircase shape. However, one can notice that at the end of every other displacement increment, the curve for $\Delta u_D = 0.1$m reaches the same value as the curves computed with the other $\Delta u_D$ values.
Table \ref{tab:crack speeds} contains the errors of the crack speeds (computed with a least-squares fit on the two numerical computations) with respect to the analytical solution.
\begin{table}[!htp]
\begin{center}
   \begin{tabular}{ | c | c | c |}
     \hline
     $\Delta u_D$ / $h$ & $0.1$ & $0.05$  \\ \hline
     $0.1$ &  $3.6\%$ & $5.1\%$ \\ \hline
     $0.01$ & $2.0\%$ & $0.68\%$  \\ \hline
     $0.001$ & $1.9\%$ & $0.70\%$ \\ \hline
   \end{tabular}
   \caption{Crack speed: error with respect to analytical solution depending on the choice of $h$ and $\Delta u_D$.}
   \label{tab:crack speeds}
\end{center}
\end{table}
The agreement of the computed crack speeds with the analytical solution is very satisfactory for all $\Delta u_D$.

\subsection{Opening mode with unknown crack path}
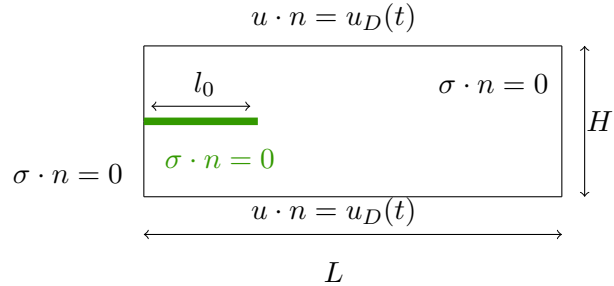
\begin{figure} [!htp]
\centering
\begin{tikzpicture}
\draw (-2.5,1) -- (3,1);
\draw (-2.5,-1) -- (3,-1);
\draw[] (-2.5,0.) -- (-2.5,1);
\draw (-2.5,0.) -- (-2.5,-1);
\draw (3,1) -- (3,-1);
\draw [<->] (-2.5,-1.5) -- (3,-1.5);
\draw (0,-2) node{$L$};
\draw[vert,line width=1mm] (-2.5,0) -- (-1,0);
\draw (-1.5,-0.5) node[vert] {$\sigma \cdot n = 0$};
\draw [<->] (3.3,-1) -- (3.3,1);
\draw (3.5,0) node{$H$};
\draw (-1.7,0.5) node{$l_0$};
\draw[<->] (-2.4,0.2) -- (-1.1,0.2);
\draw (0,1.7) node[below]{$u \cdot n = u_D(t)$};
\draw (0,-1.2) node{$u \cdot n = u_D(t)$};
\draw (2.1,0.5) node {$\sigma \cdot n = 0$};
\draw (-3.5,-0.7) node {$\sigma \cdot n = 0$};
\end{tikzpicture}
\caption{Opening mode: setup.}
\label{fig:sketch opening mode}
\end{figure}

\begin{figure}[!htp]
\centering
\includegraphics[width=0.5\textwidth]{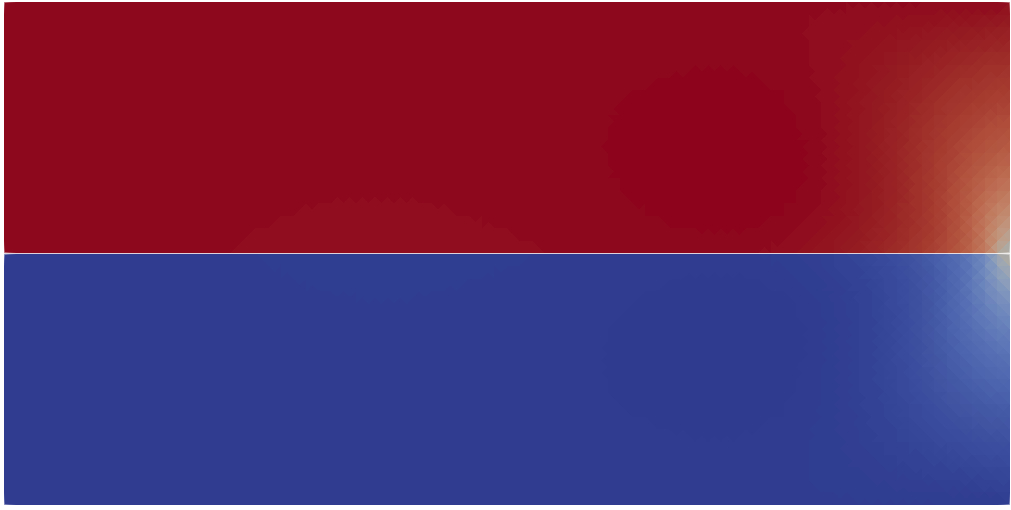}
\caption{Opening mode: $u_y$ in colors and crack path in white on a structured mesh, $u_D = 2.5$mm.}
\label{fig:mode 1 structured}
\end{figure}

The setting for this test case is presented in Figure \ref{fig:sketch opening mode}.
The dimensions of the plate are $L=32\text{mm}$ and $H=16\text{mm}$ and the initial length of the crack is $l_0 = 4\text{mm}$.
The material parameters are $E = 3.09\text{GPa}$, $\nu = 0.35$ and $\mathcal{G}_c = 300\text{kN/mm}$.
First, we use a structured mesh of size $h = 0.4\text{mm}$ leading to $25,920$ dofs. The increment in boundary conditions is defined as $\Delta u_D = h$.
Figure \ref{fig:mode 1 structured} presents the obtained crack path. We notice an unstable crack propagation, as expected, in the sense that when the propagation starts, it breaks the entire sample at a given $t_k$.

\begin{figure}[!htp]
\centering
\subfloat{
\includegraphics[width=0.5\textwidth]{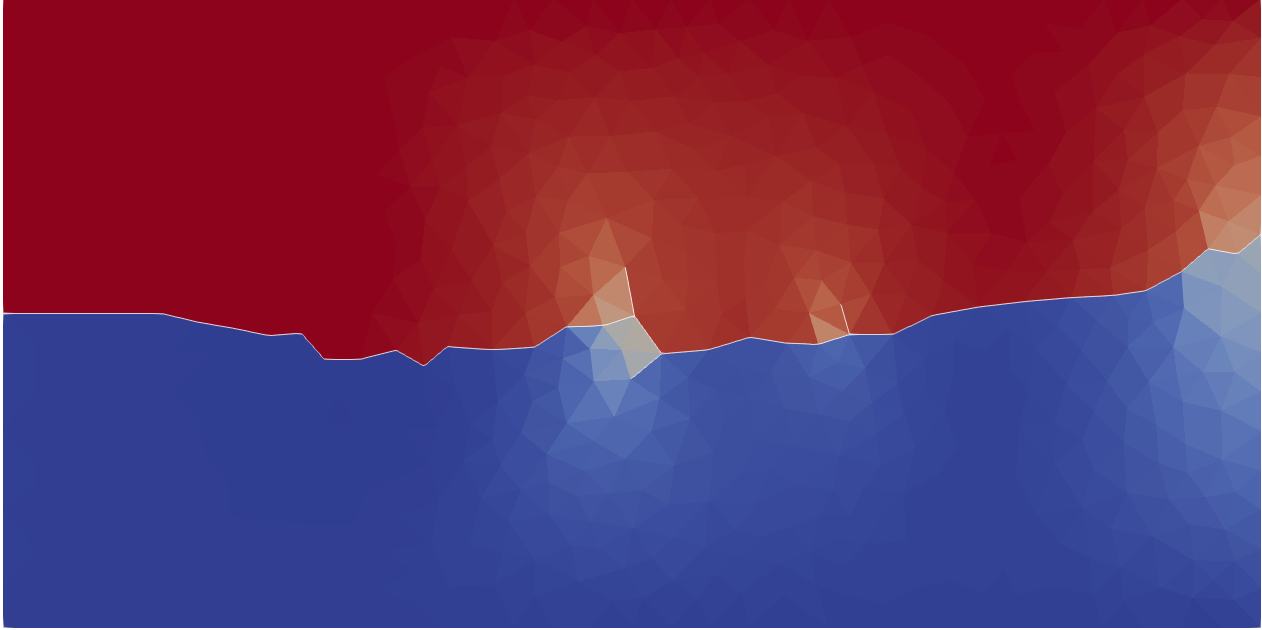}
}
\subfloat{
\includegraphics[width=0.5\textwidth]{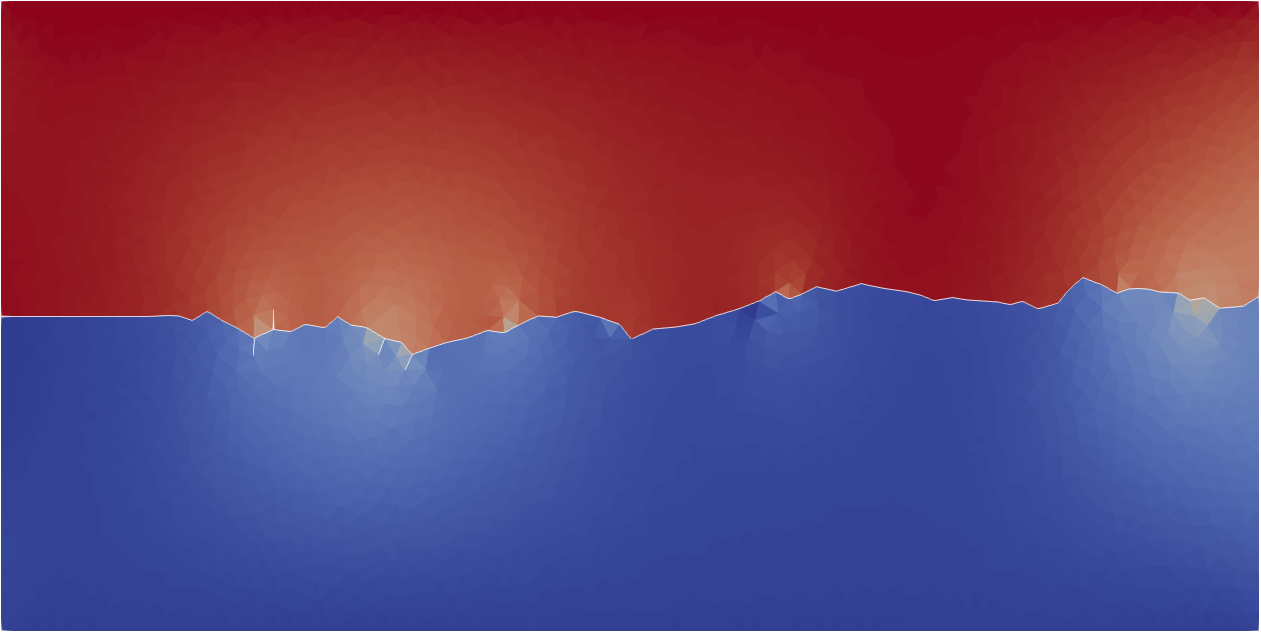}
}
\caption{Opening mode: $u_y$ in colors and crack path in white on a unstructured mesh. Left: coarse mesh, $u_D = 2.5$mm. Right: fine mesh, $u_D = 3.4$mm.}
\label{fig:mode 1 unstructured}
\end{figure}

We also perform computations on two unstructured meshes of sizes $h=1.4\mathrm{mm}$ and $h=0.74\mathrm{mm}$ corresponding respectively to $2,792$ dofs and $11,044$ dofs. Both meshes do not contain facets with a direction that could lead to a totally straight propagation of the crack.
The finer mesh is not a refinement of the coarser one.
Figure \ref{fig:mode 1 unstructured} shows the crack paths obtained on the two meshes.
The crack paths obtained are satisfactory as the propagation is rather straight and the results on the two meshes are quite similar.

\subsection{Single-edge notched shear test}
\label{sec:notched shear}
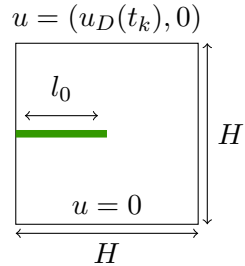
\begin{figure} [!htp]
\centering
\begin{tikzpicture} [scale=1.2]
\pgfmathsetmacro{\H}{2}
\pgfmathsetmacro{\l}{1}
\pgfmathsetmacro{\dl}{0.1}

\draw (-\H/2, -\H/2) -- (-\H/2, \H/2) -- (\H/2, \H/2) -- (\H/2, -\H/2) -- cycle;
\draw [<->] (-\H/2, -\H/2-\dl) -- (\H/2, -\H/2-\dl);
\draw[below] (0,-\H/2-\dl) node{$H$};
\draw [<->] (\H/2+\dl, -\H/2) -- (\H/2+\dl, \H/2);
\draw[right] (\H/2+\dl, 0) node{$H$};

\draw[vert,line width=1mm] (-\H/2, 0)   -- (-\H/2+\l, 0);
\draw (-0.5,0.5) node{$l_0$};
\draw[<->] (-0.9,0.2) -- (-0.1,0.2);

\draw (0,\H/2) node[above] {$u = (u_D(t_k),0)$};
\draw (0,-\H/2) node[above] {$u = 0$};
\end{tikzpicture}
\caption{Single-edge notched shear test: setup.}
\label{fig:sketch shear mode}
\end{figure}

The setting of this test case comes from \cite{ambati2015review}. It consists in a square with an already initiated crack loaded in shear on its top surface. The lower surface is recessed while the upper surface is loaded in shear. The two lateral parts are free of stress as well as the crack. Figure \ref{fig:sketch shear mode} illustrates the setting.
The crack is of initial length $l_0 = 0.5 \text{mm}$ and the dimension of the sample is $H = 1\text{mm}$.
The material parameters are $E = 210 \text{GPa}$, $\nu = 0.3$ and $\mathcal{G}_c = 2.7 \cdot 10^{-3} \text{kN/mm}$. The increment of boundary load is defined as $\Delta u_D = 10^{-6}\text{mm}$ and the final load is $u_{D,\text{final}} = 0.2 \text{mm}$. 

\begin{figure}[!htp]
\centering
\subfloat{
\includegraphics[width=0.3\textwidth]{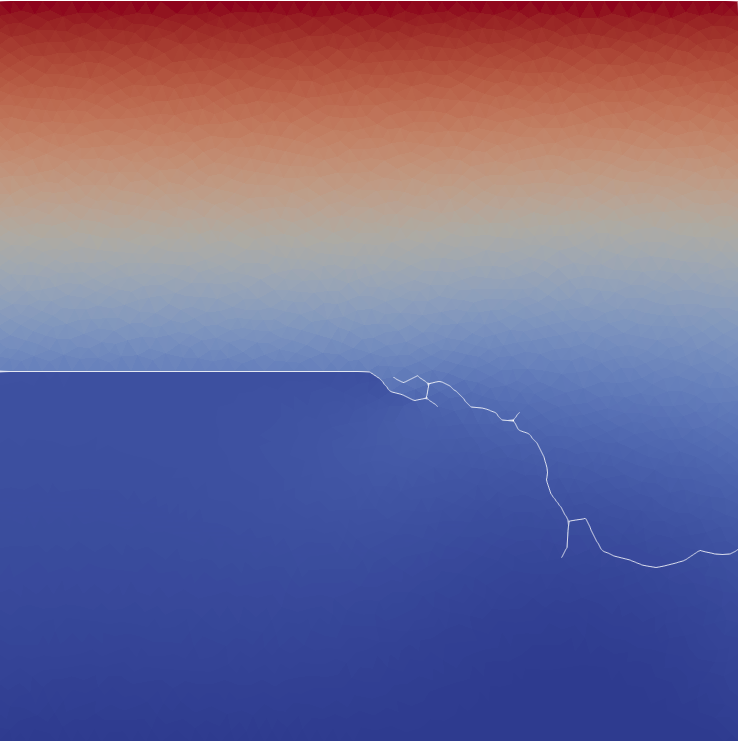}
}
\subfloat{
\includegraphics[width=0.3\textwidth]{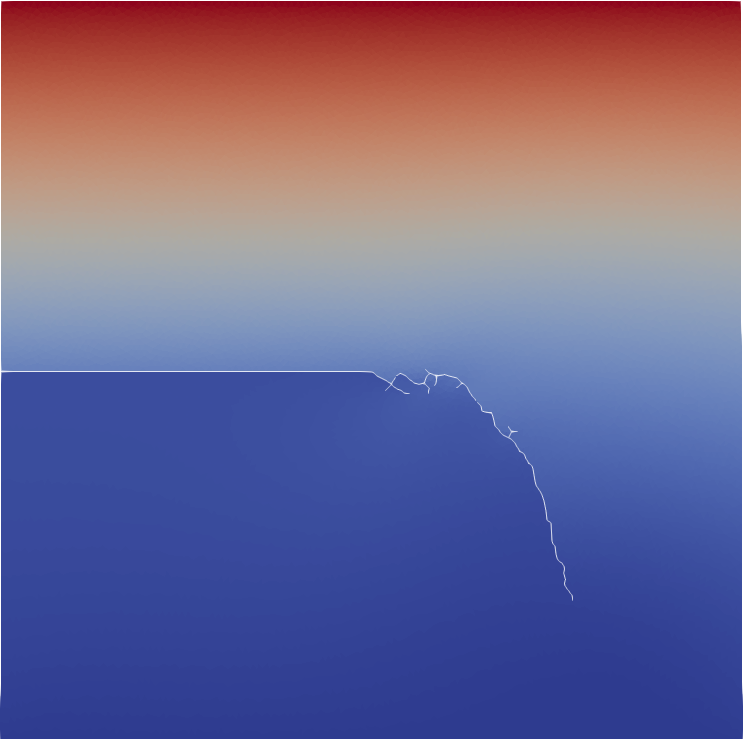}
}
\subfloat{
\includegraphics[width=0.3\textwidth]{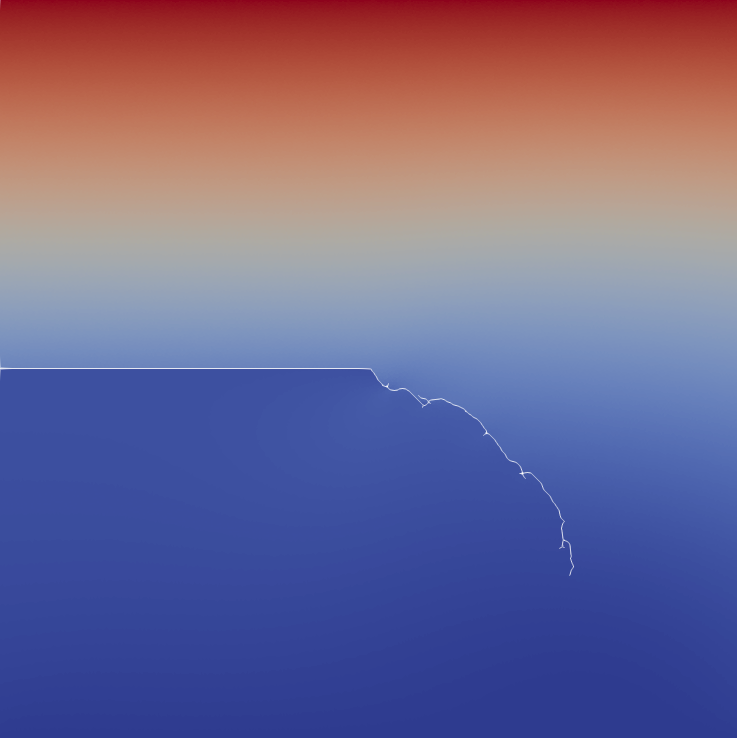}
}
\caption{Single-edge notched shear test: $u_x$ in colors and crack path in white,
$u_D = 0.2$mm. Left: coarse mesh. Middle: fine mesh. Right: finest mesh.}
\label{fig:mode 2}
\end{figure}

Three computations are performed on unstructured meshes of size $h=2.8\cdot 10^{-2}$mm (coarse mesh), $h=1.3\cdot 10^{-2}$mm (fine mesh), and $h=7.7\cdot 10^{-3}$mm (finest mesh), leading respectively to $13,396$, $65,956$, and $210,328$ dofs.
Figure \ref{fig:mode 2} shows the computed crack paths. Our results can be compared with \cite{muixi2020hybridizable} which uses a phase-field model discretized by a hybridizable discontinuous Galerkin formulation. 
The computations are in satisfactory agreement with those of \cite{muixi2020hybridizable} regarding the general orientation of the crack and the number of branches.
We observe in Figure \ref{fig:mode 2} that the crack propagates downwards along a somewhat curved path (with rather close predictions between the two finer meshes). The trajectory is sightly different from the one predicted in \cite{muixi2020hybridizable} where the crack propagates along a rather straight line which forms a sharp angle with respect to the initial crack. Experimental results would be needed to assess the correctness of these numerical results.

\begin{figure}[!htp]
\includegraphics[width=0.5\textwidth]{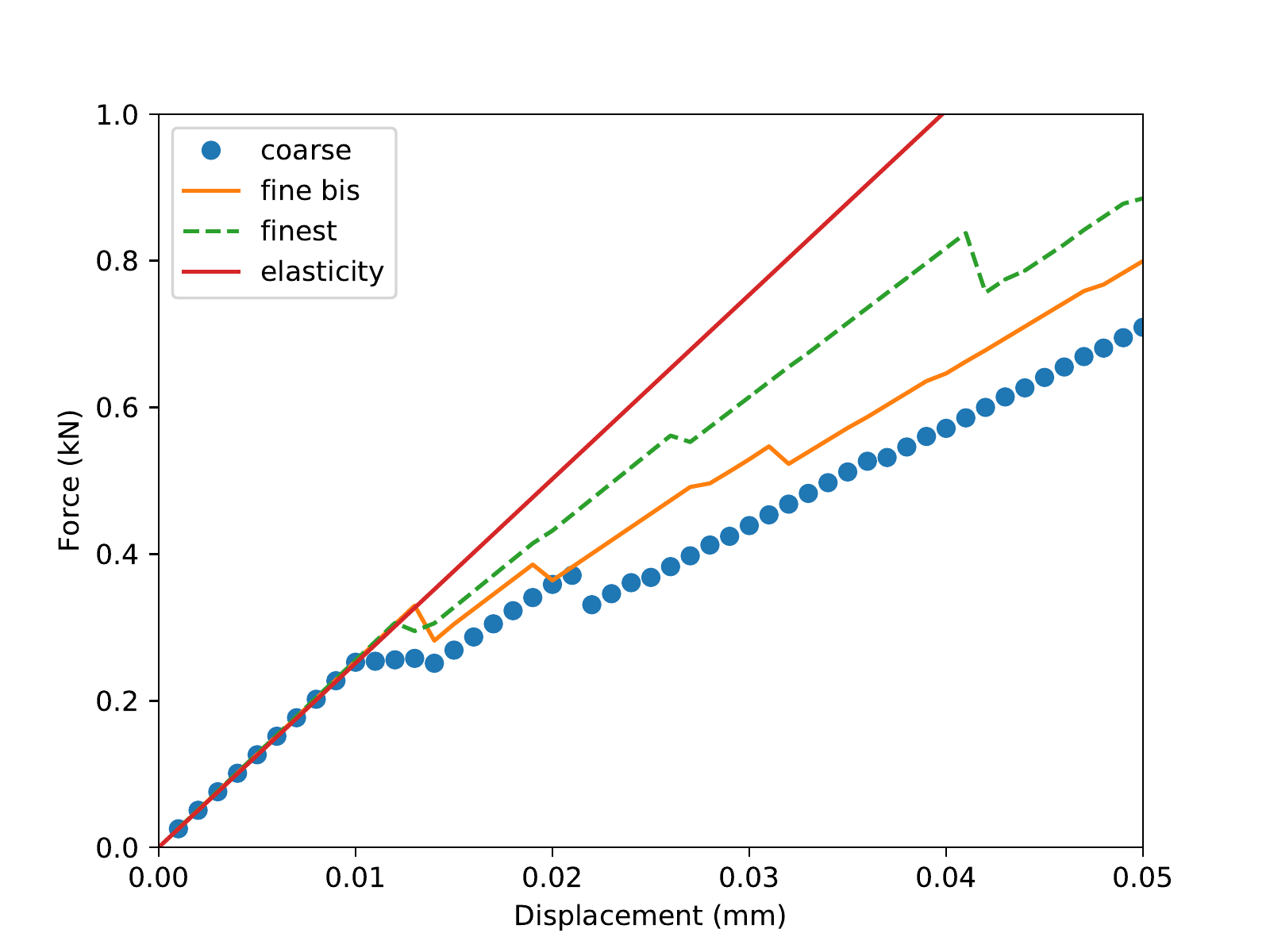}
\quad
\raisebox{3cm}{
\begin{tabular}{|c|c|c|}
\hline
mesh & $u_D$ ($\mu$m) & force (kN) \\ \hline
coarse & $9.5$ & $0.26$ \\ \hline
fine & $12.5$ & $0.33$ \\ \hline
finest & $11.3$ & $0.31$  \\ \hline
\end{tabular}}
\caption{Single-edge notched shear test. Left: load-displacement curves. Right: imposed displacement and force when the crack starts propagating.}
\label{fig:mode 2 load disp}
\end{figure}

The load-displacement curves are displayed in Figure \ref{fig:mode 2 load disp}
along with the values of the imposed displacement and the resulting force when 
the crack starts propagating.
The force is computed through an integration of the tangential component of the reconstructed normal stress $\Sigma_h \cdot n$ on the upper and lower surfaces of the sample. The force has also been computed through a residual method and the difference has been found to be negligible.
One can first notice that up to an imposed displacement of $9.5\mu$m, all the curves are superimposed and exactly reproduce the elastic response of the sample with the fixed initial crack. As the imposed displacement increases beyond the above value,
jumps in the load-displacement curves appear progressively. These jumps are a consequence of facets cracking, and the slope of the elastic response is reduced after each jump owing to the propagation of the crack. This explains the observed zigzag behavior of the response curves. Altogether, crack propagation thus induces a softening of the sample as expected. The load-displacement curve obtained on the coarse mesh stops at the value $u_D=0.2$mm for which the crack reaches the rightmost boundary of the sample. Instead, the computations on the two finer meshes support larger values for $u_D$ and lead to rather similar predictions. Furthermore,
one can see that the crack starts propagating around an imposed displacement of $10\mu$m, which is similar to the value reported in \cite{ambati2015review}. The value of the force, however, is different. We believe that this difference can be attributed to the sharp interface representation of the crack in the present method. To substantiate this claim, we performed some 
additional computations on the finest mesh using a fixed interface position, $P^1$--Lagrange finite elements, and an imposed displacement $u_D=5\mu$m. With a sharp interface, the load is $0.13$kN (consistently with the DEM prediction on the same mesh), whereas it is $0.32$kN if there is no crack (the sample is fully sound). If instead the initial crack is represented as a damage field \cite{borden2012p,muixi2020hybridizable} with a smoothing length $\ell = 5h$, the load is close to the value reported in \cite{ambati2015review,muixi2020hybridizable}, namely $0.20$kN (notice that this value is as expected in the interval $(0.13,0.32)$kN).

\subsection{Notched plate with a hole}
\begin{figure} [!htp]
\centering
\begin{tikzpicture} [scale=1.2]
\pgfmathsetmacro{\H}{4}
\pgfmathsetmacro{\L}{2}
\pgfmathsetmacro{\l}{0.5}
\pgfmathsetmacro{\dl}{0.1}
\pgfmathsetmacro{\ddl}{0.3}
\pgfmathsetmacro{\r}{0.2}
\pgfmathsetmacro{\R}{0.4}

\draw (-\L/2, -\H/2) -- (-\L/2, \H/2) -- (\L/2, \H/2) -- (\L/2, -\H/2) -- cycle;
\draw [<->] (-\L/2, -\H/2-\dl) -- (\L/2, -\H/2-\dl);
\draw[below] (0,-\H/2-\dl) node{$L$};
\draw [<->] (\L/2+\dl, -\H/2) -- (\L/2+\dl, \H/2);
\draw[right] (\L/2+\dl, 0) node{$H$};

\draw (-\L/4,\H/3) circle (\r);
\draw[<->] (-\L/4,\H/3+\ddl) -- (-\L/2,\H/3+\ddl);
\draw (-\L/4-\ddl,\H/3+1.5*\ddl) node {$a$};
\draw[<->] (-\L/4+\ddl,\H/3) -- (-\L/4+\ddl,\H/2);
\draw[right] (-\L/4+\ddl,\H/3+\ddl) node {$a$};
\draw (-\L/4,-\H/3) circle (\r);
\draw[<->] (-\L/4,-\H/3-\ddl) -- (-\L/2,-\H/3-\ddl);
\draw (-\L/4-\ddl,-\H/3-1.5*\ddl) node {$a$};
\draw[<->] (-\L/4+\ddl,-\H/3) -- (-\L/4+\ddl,-\H/2);
\draw[right] (-\L/4+\ddl,-\H/3-\ddl) node {$a$};

\draw (\L/6,0) circle (\R);
\draw[<->] (\L/6,0) -- (\L/6,\H/2);
\draw[right] (\L/6,\H/4) node {$d$};
\draw[<->] (-\L/2,0) -- (\L/6,0);
\draw[below] (-\L/4,0) node {$e$};

\draw[vert,line width=1mm] (-\L/2, \l)   -- (-\L/2+\l, \l);
\draw[<->] (-\L/2-\dl,\H/2) -- (-\L/2-\dl,\l);
\draw[left] (-\L/2-\dl, \H/4) node {$b$};
\draw (-0.7,0.9) node{$l_0$};
\draw[<->] (-1,0.7) -- (-0.5,0.7);

\end{tikzpicture}
\caption{Notched plate with a hole: setup.}
\label{fig:sketch plate holes}
\end{figure}
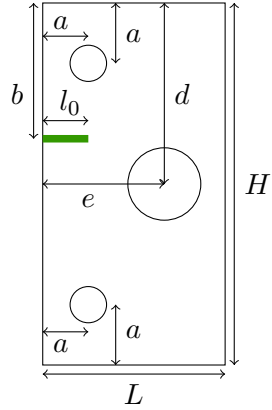

This test case comes from \cite{muixi2020hybridizable}. The material parameters are $E = 6 \text{GPa}$, $\nu = 0.22$, and $\mathcal{G}_c = 2.28\cdot 10^{-3} \text{kN/mm}$. We use fixed displacement increments of $\Delta u_D = 10^{-2} \text{mm}$.
Figure \ref{fig:sketch plate holes} presents a sketch of the sample.
The dimensions of the plate are $L=65\text{mm}$ and $H = 120 \text{mm}$.
The two holes on the left of the sample have a diameter of $10\text{mm}$ and the hole on the right of the sample has a dimeter of $20\text{mm}$. The initial length of the crack is $l_0 = 10\text{mm}$. One also has $a=20 \text{mm}$, $b=55 \text{mm}$, $d=69 \text{mm}$ and $e=36.5 \text{mm}$.
The right hole is free of stress, the lower hole is recessed and the upper hole has an imposed displacement $u = (0,u_D(t_k))$.
We use three unstructured meshes with $h=2.8$mm, $h=1.5$mm and $h=0.78$mm having respectively $9,926$, $39,380$ and $157,340$ dofs.

\begin{figure}[!htp]
\centering
\subfloat{
\includegraphics[width=0.3\textwidth]{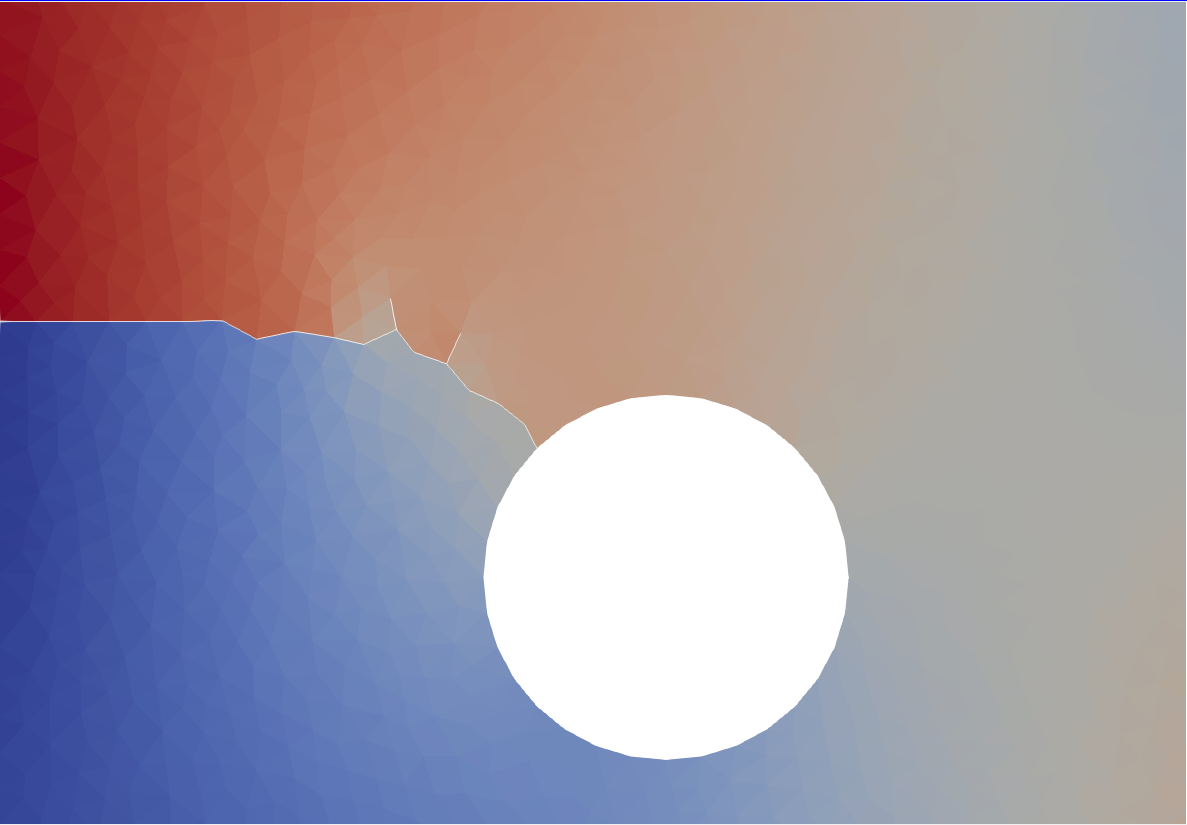}
}
\subfloat{
\includegraphics[width=0.3\textwidth]{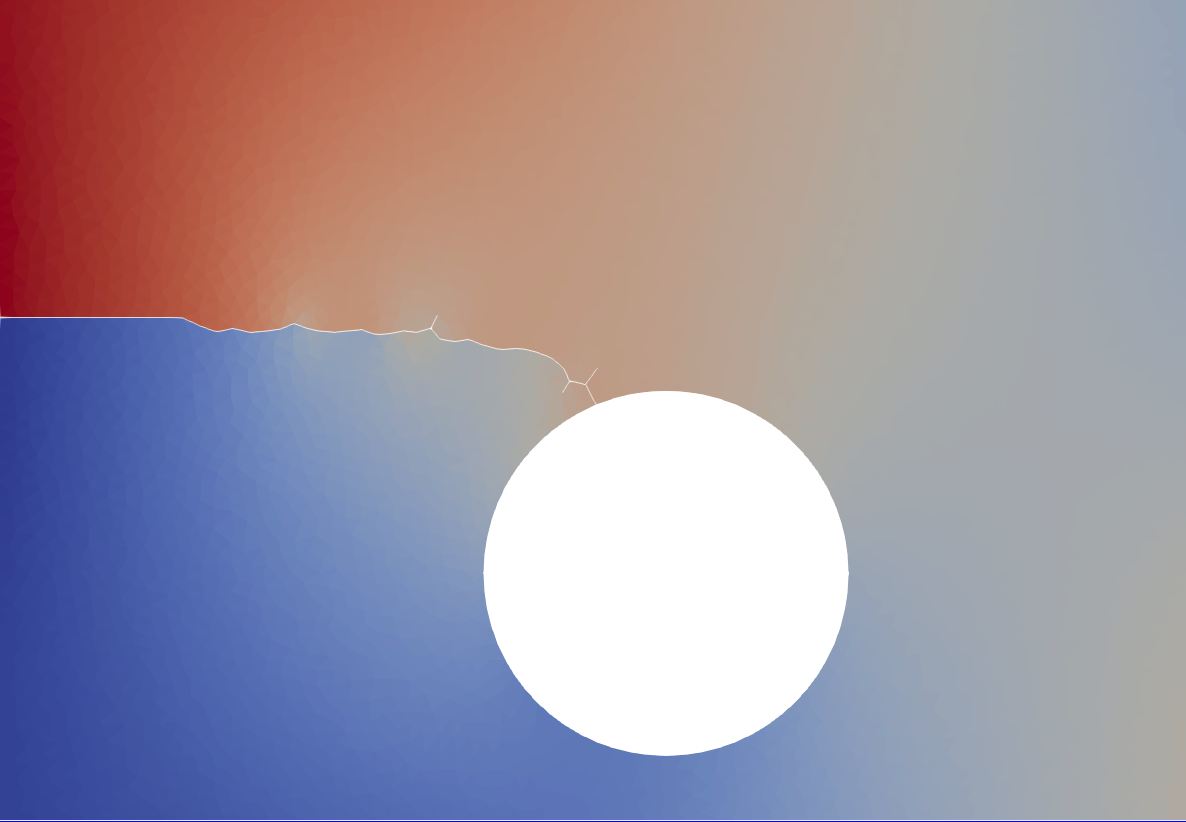}
}
\subfloat{
\includegraphics[width=0.3\textwidth]{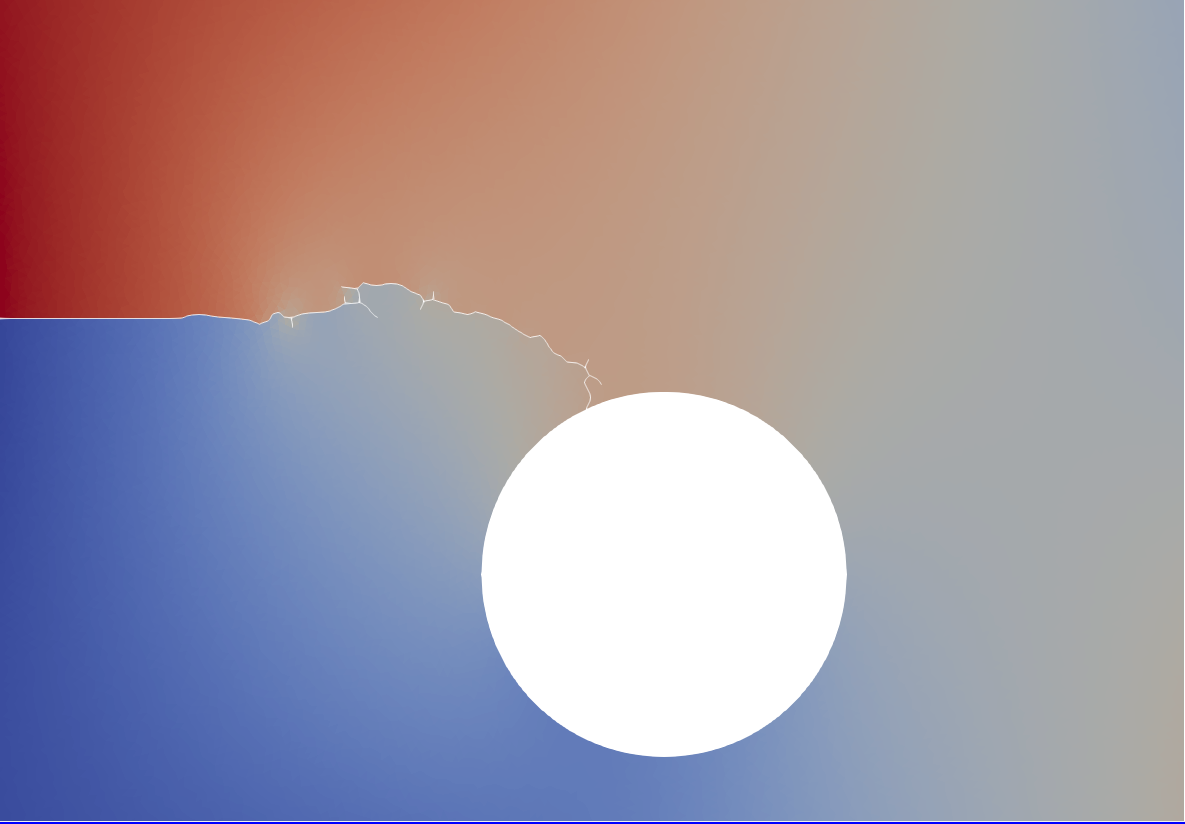}
}
\caption{Notched plate with a hole: $u_y$ in colors and crack path in white (zoom around the main hole). Left: coarse mesh, $u_D=0.65$mm. Middle: fine mesh, $u_D=1.4$mm. Right: finest mesh, $u_D=3.0$mm.}
\label{fig:plate holes}
\end{figure}

\begin{figure}[!htp]
\centering
\includegraphics[width=0.5\textwidth]{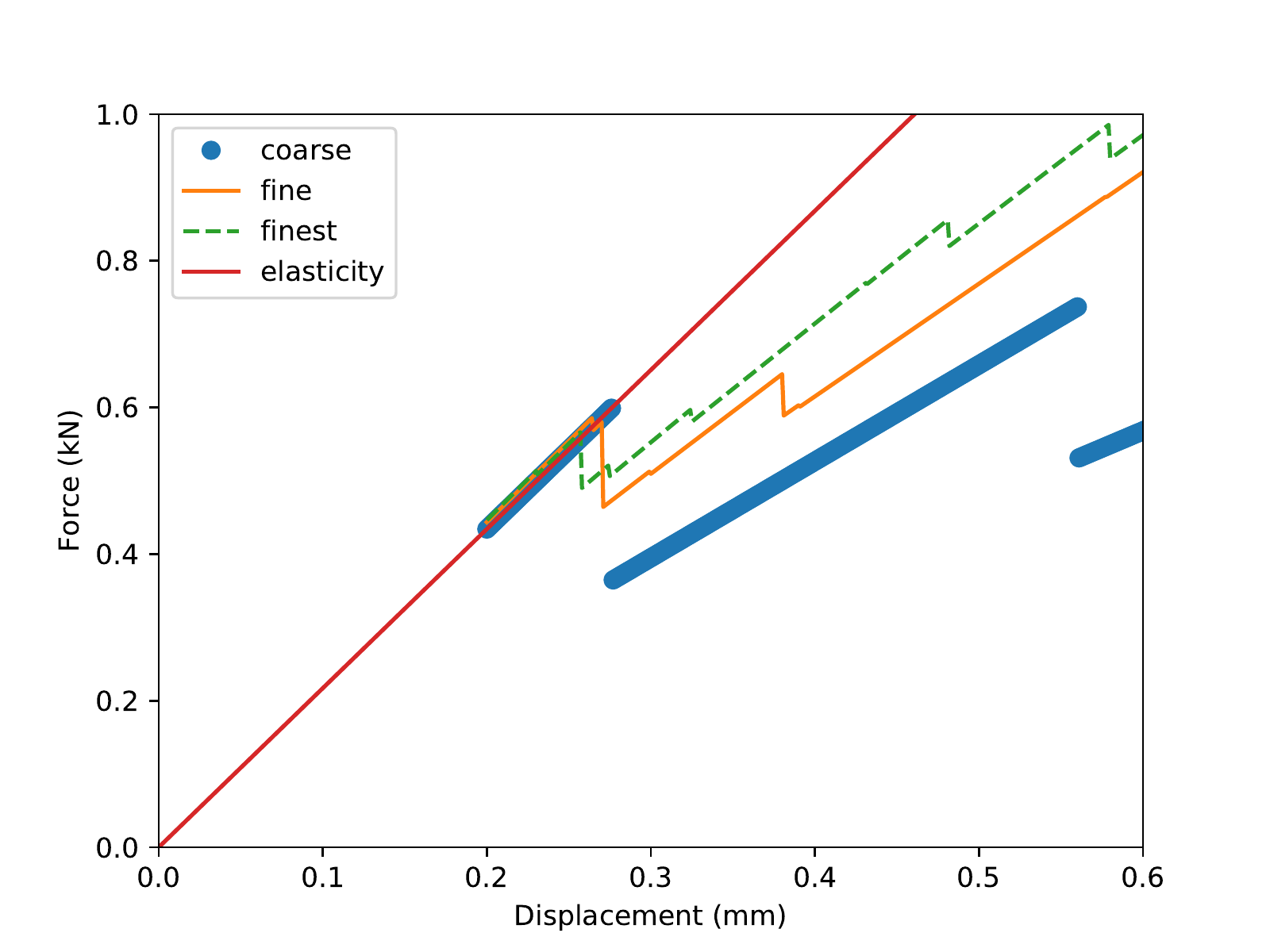}
\quad
\raisebox{3cm}{
\begin{tabular}{|c|c|c|}
\hline
mesh & $u_D$ (mm) & force (kN) \\ \hline
coarse & $0.28$ & $0.60$ \\ \hline
fine & $0.27$ & $0.58$ \\ \hline
finest & $0.26$ & $0.57$ \\ \hline
\end{tabular}}
\caption{Notched plate with a hole. Left: load-displacement curves. Right: imposed displacement and force when the crack starts propagating.}
\label{fig:plate holes load disp}
\end{figure}

\begin{figure}[!htp]
\centering
\includegraphics[width=0.5\textwidth]{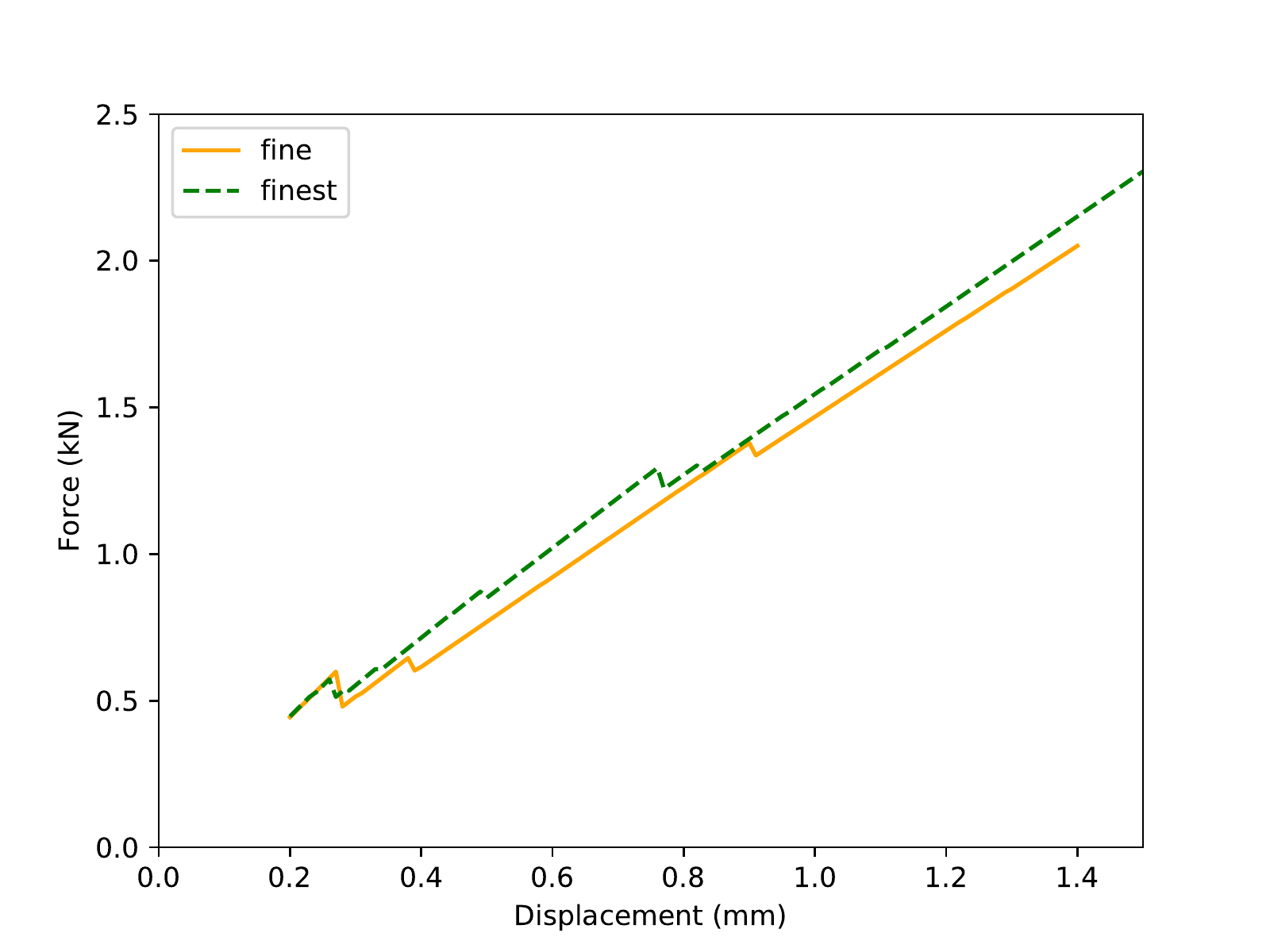}
\caption{Notched plate with a hole: load-displacement curves.} 
\label{fig:plate holes load disp 2}
\end{figure}

Figure \ref{fig:plate holes} shows the computed crack paths.
We compare our results with \cite{muixi2020hybridizable} without taking into account the secondary crack starting from the largest hole as we restrict ourselves to crack propagation and not crack initiation.
We notice that for the three computations, the crack goes towards the largest hole in a similar fashion which also seems consistent with \cite{muixi2020hybridizable}.
The load-displacement curves are given in Figures \ref{fig:plate holes load disp} and \ref{fig:plate holes load disp 2}, together with the values of the imposed displacement and the resulting force when the crack starts propagating and when the crack reaches the hole, respectively. Figure~\ref{fig:plate holes load disp} focuses on imposed displacements $u_D$ up to 0.6mm, whereas Figure~\ref{fig:plate holes load disp 2} explores a wider range for $u_D$ on the two finer meshes.
The force is computed through an integration of the vertical component of the reconstructed normal stress $\Sigma_h \cdot n$ on the upper left hole.
A similar behaviour of the elastic response and the softening of the sample is observed as in Section \ref{sec:notched shear}. The crack starts propagating around an imposed displacement $u_D = 0.27$mm (consistently on the three meshes), which is in reasonable agreement with the caption of \cite[Fig.~19]{muixi2020hybridizable} which indicates that propagation has started at the value of $u_D=0.3$mm.
A further quantitative comparison including forces is delicate owing to the difficulties mentioned at the end of Section \ref{sec:notched shear}.
Moreover, we observe from Figure~\ref{fig:plate holes load disp} that the predictions on the coarse mesh are still rather inaccurate for higher values of $u_D$, whereas Figure~\ref{fig:plate holes load disp 2} indicates that the predictions on the two finer meshes are in satisfactory agreement as far as the load-displacement curves are concerned.
The predictions of the path of crack propagation are also similar on both meshes, but the value of the imposed displacement when the crack reaches the hole is different, as reflected in the caption of Figure~\ref{fig:plate holes}.

\section{Conclusions}
\label{sec:conclusions cracking}
We have presented a variational Discrete Element Method (DEM) to compute Griffith crack propagation. The crack propagates through the facets of the mesh and thus between discrete elements.
The variational DEM is a consistent discretization of a Cauchy continuum and only requires three continuum macroscopic parameters for its implementation: the Young modulus, the Poisson ratio, and the critical energy release rate. The displacement degrees of freedom are attached to the barycentre of the mesh cells.
A discrete Stokes formula is used to devise a piecewise constant gradient and linearized strain reconstructions.
An approximation of the energy release rate is computed in the procedure \texttt{ESTIMATE}. The procedure \texttt{MARK} then determines the breaking facet at each pseudo-time node $t_k$. Finally, the procedure \texttt{UPDATE} updates the necessary discrete quantities after the facet that has been marked has been broken. A convergence test in antiplane shear has confirmed the efficiency of the variational DEM discretization as well as the $\mathcal{O}(h^{\frac12})$ convergence rate in energy norm.
The robustness of the method regarding the computation of the crack speed has been verified. Also, several numerical experiments have shown that the method can provide reasonable crack paths.

This work can be pursued in several directions. A first idea would be to adapt the present methodology to three-dimensional problems with two-dimensional cracks.
A second direction concerns the regularity of the crack surface. Indeed, in the spirit of \cite{francfort1998revisiting}, a crack should be a surface that minimizes energy. To achieve this goal, the variational DEM could be coupled to gradient flows used for surface lifting, as in \cite{romon2013introduction}, with the goal of moving the crack surface vertices. One would then have to verify the convergence of the discrete crack area with tools similar to \cite{hildebrandt2006convergence}.
A third direction for further study is to approximate cohesive cracking laws instead of a Griffith cracking law so as to enable the simulation of crack initiation as well as crack propagation. Inspiration can be found in \cite{mariotti2009modeling} which uses a DEM with a linear cohesive law.
Finally, a last direction can be to consider an enrichment similar to \cite{chahine2008crack} close to the crack tip so as to obtain a convergence with order $\mathcal{O}(h)$.

\section*{Acknowledgements}
Partial support by CEA is gratefully acknowledged.

\bibliographystyle{plain} 
\bibliography{bib_globale} 

\end{document}